\documentstyle[10pt]{amsart}
\newcommand{\R}{{\Bbb R}}
\newcommand{\Z}{{\Bbb Z}}
\newcommand{\Nb}{{\Bbb N}}

\newcommand{\C}{{\Bbb C}}

\newcommand{\half}{\frac{1}{2}}

\theoremstyle{plain}

\title{Wave invariants for non-degenerate closed geodesics}
\author{Steven Zelditch*}
\thanks{ *Partially supported by NSF grant \#DMS-9404637.}
\address{Johns Hopkins University, Baltimore, Maryland  21218 \newline
email:{\tt zel$\mbox{\@}$}chow.mat.jhu.edu}  
\begin{document}
\maketitle

\addtolength{\baselineskip}{1pt} 

\setcounter{section}{-1}
\section{Introduction}

This paper is a continuation of [Z.1]. There,
we gave an effective method for putting the
Laplacian $\Delta$ of a Riemannian manifold $(M,g)$ into a quantum Birkhoff
normal form around a non-degenerate elliptic closed geodesic
$\gamma$, and applied it to the calculation and characterization of the 
wave invariants $a_{\gamma k}$ at $\gamma.$  The wave invariants, we recall, are
 the coefficients in the
singularity expansion
\medskip

\begin{tabular}{l}$Tr U(t) = e_o(t) + \sum_{\gamma}  e_{\gamma}(t)$\\
$e_{\gamma}(t) \sim a_{\gamma\;-1}(t - L_{\gamma} + i0)^{-1} +
\sum_{k=0}^{\infty} a_{\gamma\;k}(t-L_{\gamma} +i0)^k log (t - L_{\gamma} +
i0)$
\end{tabular} 
\medskip

\noindent of the trace of the wave group $U(t) =
e^{it\sqrt{\Delta}}$ at lengths $t = L_{\gamma}$ of closed geodesics. 
 The first purpose of this article is to
extend the methods and results of [Z.1] to  general non-degenerate closed geodesics,
i.e. to $\gamma$ whose Poincare map $P_{\gamma}$ is any symplectic sum of 
non-degenerate elliptic,  hyperbolic, or loxodromic  parts. Our second purpose is to 
generalize to the full non-degenerate case
 the inverse result of Guillemin that the quantum normal form coefficients at
non-degenerate elliptic  closed geodesics are spectral invariants [G.1,2].  
 It will 
follow that, for metrics with simple length spectra, the classical Birkhoff normal
form of the metric around any non-degenerate closed geodesic is a spectral invariant
of the Laplacian.

Let us state the results  precisely.  The first is that   
the wave invariants $a_{\gamma k}$ in the general non-degenerate case are essentially
  analytic continuations of the expressions obtained in the elliptic case.
They may be written in the form 
$$a_{\gamma k} =  {\cal F}_{k, -1}(D)\cdot Ch (x)|_{x = P_{\gamma}}\leqno(0.1)$$
where 
$$Ch (x) = \frac{i^{\sigma}}{\sqrt{|det (I - x)|}}$$
is the character 
of the metaplectic representation (with $\sigma$  a certain
Maslov index) and where ${\cal F}_{k, -1}(D) $ is an invariant 
 partial differential operator on the metaplectic group $Mp(2n, \R)$ which
is canonically fashioned from the germ of the metric $g$ at $\gamma$.  The
exact expression for ${\cal F}_{k, -1}(D)$ will be given in \S 5 and leads to
the following characterization of the wave invariants (cf. [Z.1, Theorem A]):
 \medskip

\noindent{\bf Theorem I} {\it  Let ${\gamma}$ be a non-degenerate closed
geodesic.  Then 
$a_{\gamma k} = \int_{\gamma} I_{\gamma; k} (s; g)ds$ 
where:

\noindent(i)  $I_{\gamma; k}(s,g)$  is a homogeneous Fermi-Jacobi-Floquet
 polynomial  of weight -k-1 in the data
$\{ y_{ij}, \dot{y}_{ij}, D^{\beta}_{s,y}g \}$ with $|\beta| \leq 2k+4$ ;

\noindent(ii) The degree of $I_{\gamma; k }$ in the 
Jacobi field components is
at most 6k+6;

\noindent(iii) At most 2k+1 
indefinite integrations over $\gamma$ occur in $I_{\gamma k r}$;

\noindent(iv) The degree of $I_{\gamma; k }$ 
in the Floquet invariants $\beta_j$ is
at most k+2.}
\medskip

 The relevant terminology and notation will be recollected in  \S 1.  From (0.1)
and from the formulae for ${\cal F}_{k, -1}(D)$ and $Ch$, we  give a rather
simple proof (\S 6) of the following inverse result (strictly speaking, proven
 only for non-degenerate elliptic closed geodesics in [G.2, Theorem 1.4]):
\medskip

\noindent{\bf Theorem II } {\it Let $\gamma$ be a
non-degenerate closed geodesic. Then
the entire quantum Birkhoff normal form around $\gamma$ is a spectral invariant;
in particular the classical Birhoff normal form is a spectral invariant.}
\medskip

We thus have:
\medskip

\noindent{\bf Corollary II.1} {\it Suppose $(M,g)$ is a compact Riemannian
manifold with simple length spectrum $Lsp(M,g)$ and with all closed geodesics
non-degenerate.
  Then from
$Spec(M,g)$ one can recover the  quantum (and hence classical) Birkhoff normal
forms around all closed geodesics.}
\medskip

The hypotheses of the corollary are of course satisfied by  generic Riemannian
metrics (cf. [Kl, Lemma 4.4.3]).   The corollary therefore answers affirmatively
the third question in [Z.2, p.692], which asks whether isospectral manifolds in
this class of metrics are locally Fourier isospectral near corresponding closed
geodesics.  

Let us now briefly discuss the main ideas in the proofs, and in particular the novel
aspects caused by the hyperbolic and loxodromic parts of $P_{\gamma}$.

As in [G.1,2][Z.1], the wave invariants at a closed geodesic $\gamma$ will be
expressed as non-commutative residues of the wave group and its time derivatives
at $t = L_{\gamma}.$ For simplicity we will often abbreviate $L_{\gamma}$ by $L$.
Then we have:
 $$a_{\gamma k} = res\ D_t^k e^{it\sqrt{\Delta}}|_{t = L}: = Res_{s=0} TrD_t^k
 e^{it\sqrt{\Delta}} \sqrt{\Delta}^{-s}|_{t = L}.$$ 
Since $res$  is invariant under conjugation by (microlocal) unitary Fourier integral
operators, the $a_{\gamma k}$'s may be
calculated by putting the wave group into a microlocal (quantum Birkhoff) normal form
around $\gamma$ and by 
 by determining the residues of the resulting wave group of
the normal form.

The primary step in the  analysis of the wave invariants $a_{\gamma
k}$ is therefore to put $\Delta$ into this microlocal normal form around
$\gamma.$ In the case of non-degenerate elliptic closed geodesics, we recall, the
normal form was a polyhomogeneous function in the (microlocally elliptic) element
$D_s$ with coefficients in the transverse (elliptic) harmonic
 oscillators $\hat{I}^e_j = \half (D_{y_j}^2 + y_j^2)$ [Z.1,Theorem B].  
  In the general non-degenerate case, the normal
form will involve a greater variety of quadratic normal
forms or `action operators:' in addition to the elliptic action operator $\hat{I}^e_j$
there  can also occur the real hyperbolic action operators $\hat{I}^h_j$ and
complex hyperbolic (or loxodromic) action operators $\hat{I}_j^{ch, Re},
\hat{I}_j^{ch, Im}.$  These hyperbolic actions cause several complications to
the arguments in the elliptic case: First, they have continuous spectra, and
so the construction of the intertwining operator to the normal form has to
be modified in several ways (\S 3,4).  Second, the wave group of the normal form has
continuous spectrum and this alters the calculation of its residues (\S 5).  Third,
the presence of real parts in the Floquet exponents of $P_{\gamma}$ complicates
the process of determining the normal form coefficients from the wave invariants
(\S 6).

To get acquainted with these action operators and the normal form algorithm, let us
consider the very first step of  ``linearizing" 
$\sqrt{\Delta}$  around $\gamma$ and of putting the `linearization'
$${\cal L} = D_s - \half (\sum_{j=1}^n D_{y_j}^2 + \sum_{ij=1}^n K_{ij}y_iy_j)$$
into quantum quadratic normal form.  Here,  $ n = dim M -1$, the coordinates $(s,y_j)$
are the (re-scaled) Fermi normal coordinates around $\gamma$, $D_s =
\frac{\partial}{i \partial s} $ and
$K_{ij}$ is the curvature operator $g(R(\partial_s, \partial_{y_i})\partial_s,
\partial_{y_j})$.  The linearization ${\cal L}$ of $\sqrt{\Delta}$ is a quadratic
Hamiltonian and is the Weyl quantization of a quadratic classical Hamiltonian
(see [Ho III] and [Ho] for background on Weyl quantizations and normal forms for
quadratic Hamiltonians).  Hence its symbol may be conjugated into normal form 
by an element of ${\cal W} \in Sp(2n, \R).$  The operator ${\cal L}$ itself 
 may be put into normal form by 
conjugating with the metaplectic operator $\mu({\cal W})$ quantizing  ${\cal W}$.  As
will be seen in \S1 , this linear symplectic map  is   the Wronskian matrix ${\cal W}$
whose columns consist of the Jacobi eigenfields of $P_{\gamma}$.
The normal form of the linearized $\sqrt{\Delta}$ is therefore given by ${\cal R} =
\mu({\cal W})^* {\cal L} \mu({\cal W})^{* -1}.$

 In the elliptic
case [Z.1] $P_{\gamma}$  was a direct sum of rotations, and the quantum normal form
of ${\cal L}$ had the form
$${\cal R}^e = D_s + \frac{1}{L} H_{\alpha},\;\;\;\;\;\;\; H_{\alpha} =
\sum_{j=1}^n \alpha_j \hat{I}_j^e$$
where the spectrum $\sigma(P_{\gamma}) = \{e^{\pm i \alpha_j}\}$. In the general
non-degenerate case the normal form will similarly depend on   the spectral
decomposition of
$P_{\gamma}$. Recall that,  
 since $P_{\gamma}$ is 
symplectic, its eigenvalues $\rho_j$  come in three types: (i) pairs 
$\rho, \bar{\rho}$ of conjugate eigenvalues of modulus 1; (ii) pairs $\rho,
\rho^{-1}$ of inverse real eigenvalues; and (iii) 4-tuplets $\rho, \bar{\rho},
\rho^{-1} \bar{\rho}^{-1}$ of
complex eigenvalues. We will often write them in the forms: (i) $e^{ \pm i\alpha_j} $,
(ii)$e^{\pm \lambda_j}$, (iii)  $e^{\pm  \mu_j \pm i \nu_j}$ respectively
 (with $\alpha_j, \lambda_j, \mu_j, \nu_j \in \R$), although
a pair of inverse real eigenvalues $\{-e^{\pm \lambda}\}$ could  be negative.  Here,
and throughout, we make the assumption that $P_{\gamma}$ is {\it non-degenerate}
in the sense that
$$\Pi_{i=1}^{2n} \rho_i^{m_i} \not= 1,\;\;\;\;\;\;\;\;\;(\forall \rho_i \in
\sigma(P_{\gamma}),\;\;\;\;\; (m_1,\dots,m_{2n}) \in {\bf N}^{2n}).$$
Each type of eigenvalue then determines a different type of quadratic
action, both on the classical and quantum levels (cf. [Ho, Theorem 3.1],[Ar]):
\medskip

\begin{tabular}{l|l|l} Eigenvalue type & Classical Normal form & Quantum normal form
\\
\hline   (i) Elliptic type  & $I^e = \half \alpha (\eta^2 + y^2)$
&  
$\hat{I}^e:=
\half \alpha (D_y^2 +y^2)$
\\
$\{e^{i\pm \alpha}\}$ & & \\
\hline (ii) Real hyperbolic type & $I^{h} = 2 \lambda y \eta$ & $\hat{I}^{h}:= 
\lambda (y D_y + D_y y)$
\\
$\{e^{\pm \lambda}\}$ & & \\ 
\hline(iii) Complex hyperbolic   & $I^{ch,Re} = 2 \mu (y_1 \eta_1 + y_2 \eta_2)$ &
 $\hat{I}^{ch, Re} = \mu  (y_1D_{y_1} + D_{y_1} y_1 + y_2 D_{y_2} + D_{y_2} y_2),$\\
(or loxodromic type) & $I^{ch,Im} = \nu ( y_1 \eta_2 - y_2 \eta_1)$ & $\hat{I}^{ch,Im}
=   \nu (y_1D_{y_2}  - y_2 D_{y_1} )$
\\
$\{e^{\pm \mu + i\pm \nu}\}$ & & \\
\end{tabular}
\medskip 

  In the case where  the Poincare map $P_{\gamma}$ has p
pairs of complex conjugate eigenvalues of moduls 1, q pairs of inverse real
eigenvalues and c quadruplets of complex hyperbolic eigenvalues,
the linearized $\sqrt{\Delta}$ will have the form:
$${\cal R} =  D_s +  \frac{1}{L}[\sum_{j=1}^p \alpha_j
\hat{I}_j^e + \sum_{j=1}^q \lambda_j \hat{I}_j^h + \sum_{j=1}^c \mu_j \hat{I}^{ch,
Re}_j +
\nu_j \hat{I}_j^{ch, Im}]. $$
The full
quantum Birkhoff normal form is then given by the analogue of Theorem B
of [Z.1]:
\medskip

\noindent{\bf Theorem B} ~~~~{\it Assuming $\gamma$ non-degenerate, there exists a
microlocally elliptic Fourier integral operator $W$ from the conic neighborhood of
$\R^+ \gamma$ in $T^*(N_{\gamma})$ to  the corresponding cone in $T_+^*S^1$ in
$T^*(S^1
\times
\R^n)$ such that
$$W \sqrt{\Delta} W^{-1} \equiv D_s + \frac{1}{L}[\sum_{j=1}^p \alpha_j
\hat{I}_j^e + \sum_{j=1}^q \lambda_j \hat{I}_j^h + \sum_{j=1}^c \mu_j \hat{I}^{ch,
Re}_j +
\nu_j \hat{I}_j^{ch, Im}] +$$
$$ +  \frac{p_1(\hat{I}_1^e,\dots,
\hat{I}_p^e, \hat{I}_1^h, \dots, \hat{I}_q^h, \hat{I}_1^{ch,
Re},\hat{I}_1^{ch,Im},\dots, \hat{I}_c^{ch,Re}, \hat{I}_c^{ch, Im})}{D_s} +\dots$$
$$+ \frac{p_{k+1}(\hat{I}_1^e,\dots, \hat{I}_c^{ch, Im})}{D_s^k} +\dots $$
where the numerators $p_j(\hat{I}_1^e,\dots, \hat{I}_p^e, \hat{I}_1^h, \dots,
\hat{I}_c^{ch, Im})$ are polynomials of degree j+1 in the variables
$(\hat{I}_1^e,\dots,  \hat{I}_c^{ch, Im})$ and where the kth remainder term lies in
the space $\oplus_{j=o}^{k+2} O_{2(k+2-j)}\Psi^{1-j}$}        
\medskip

Here, $O_n \Psi^r$ is the space of pseudodifferential operators of order
r whose complete symbols vanish to order n at $(y,\eta)=(0,0).$  Thus, the
remainder terms are `small' in that they combine in some mixture  a low
pseudodifferential order or a high vanishing  order along $\gamma$. 

Some remarks now on the contents and organization of this paper. Since
it is a continuation of [Z.1], we have tried to avoid duplicating arguments and
calculations which are essentially unchanged from the elliptic case. 
Many of the arguments which remain are still quite analogous to the
elliptic case and inevitably produce a sense of deja-vu. Our excuse for
drawing them out to their present length is that it is not apriori clear that
the arguments of the elliptic case  generalize so neatly to the hyperbolic
and loxodromic cases.   It may even be viewed as a virtue of the method of [Z.1] that
it adapts so effortlessly to the general case.

 It should be noted here that Guillemin was aware that
the arguments of the elliptic case should extend to the general non-degenerate
case and stated his main result, Theorem 1.4 of [G.2], for the general case.
However, we also note that the methods used here are  extensions of the
methods of [Z.1], which in many significant respects differ from the methods of
[G.2].

The organization of this paper is as follows: In \S 1 we will review the symplectic
and microlocal ingredients required to construct a `linear model' for the Laplacian
near a closed geodesic $\gamma$. In \S 2, we introduce the semi-classically scaled
Laplacian and the linearized Laplacian and conjugate them to the model space.  In \S
3, we conjugate the resulting  semi-classical model Laplacian to a semi-classical
normal form to infinite order.  In \S 4 we show how this semi-classical normal form
induces a bona-fide quantum Birkhoff normal form for the Laplacian near $\gamma$.
In \S 5, we use the normal form to give the  explicit formula (0.1) for the wave
invariants. In \S 6, we show that the
 quantum normal form coefficients can be determined  from the special
values of (0.1) corresponding to $\gamma$
and its iterates.

\section{Preliminaries}

This section begins with a resume of the symplectic linear algebra underlying the
Jacobi equation, the linear Poincare map, and the symplectic classification of
non-degenerate quadratic forms (\S1.1).  It then summarizes the quantum aspects
of the linear theory, in particular the behaviour of the quantum action operators
(\S 1.2)

\noindent{\bf \S1.1: Symplectic preliminaries}\\

\noindent{\bf \S1.1a: Closed geodesics, linear Poincare maps and Jacobi fields}\\

Throughout this paper, $\gamma$ will  denote 
 a {\it primitive} closed geodesic of $(M,g)$; its iterates will be denoted by
$\gamma^m$. 

The space   ${\cal J}_{\gamma}^{\bot}$ of (real) orthogonal Jacobi fields
along
$\gamma$ is then the real symplectic vector space, of  dimension 2n, of solutions
of the Jacobi equation $Y'' + R(T,Y)T = 0$ (with $T$ the unit tangent vector
along $\gamma.$)  The symplectic structure is given by the Wronskian
$$\omega(X,Y) = g(X, \frac{D}{ds}Y) - g(\frac{D}{ds}X, Y).$$

 The linear Poincare map $P_{\gamma}$ is  the (real) linear symplectic map 
on $({\cal J}_{\gamma}^{\bot},\omega)$ defined by $P_{\gamma} Y(t) = Y(t +
L_{\gamma}).$  To diagonalize it, we also  
 complexify it as a complex symplectic map $P_{\gamma}^{\C}$ on
the space 
 ${\cal J}_{\gamma}^{\bot}\otimes \C$ of complex orthogonal Jacobi fields. 
Here, the symplectic form is extended to the complexified space as a complex bilinear
form
$\omega^{\C}$.
 Since $P_{\gamma}^{\C} \in Sp({\cal J}_{\gamma}^{\bot}\otimes \C,\omega)$,
its spectrum $\sigma (P_{\gamma}^{\C})$ is stable under inverse and complex
conjugation: thus, if $\rho \in \sigma (P_{\gamma}^{\C})$, then also $ \rho^{-1},
\bar{\rho},
\bar{\rho}^{-1} \in \sigma (P_{\gamma}^{\C})$. As mentioned above, we will assume that
$P_{\gamma}^{\C}$ is {\it non-degenerate} in the following  strong sense:
$$\rho_1^{m_1} \dots \rho_n^{m_n} = 1 \Rightarrow m_i = 0\;\;\;\; (\forall i, m_i \in 
\Nb ).\leqno(1.1a.1)$$
In particular, the eigenvalues are simple and $\pm 1 \notin \sigma
(P_{\gamma}^{\C}).$

The eigenspace of $P_{\gamma}^{\C}$ of eigenvalue $\rho$ will be denoted by 
${\cal J}_{\gamma}^{\bot, \C}(\rho) \subset {\cal J}_{\gamma}^{\bot}\otimes \C.$
We then have the symplectic orthogonal decomposition 
$${\cal J}_{\gamma}^{\bot, \C}\otimes \C = {\cal J}_{nc}^{\bot, \C} \oplus {\cal
J}_{co}^{\bot, \C}\leqno(1.1a.2)$$
into the `non-compact'  symplectic subspace  
$${\cal J}_{nc}^{\bot, \C} = \oplus_{\rho: |\rho|\not= 1} {\cal
J}_{\gamma}^{\bot, \C}(\rho) \leqno(1.1a.3.nc)$$ 
where $P_{\gamma}^{\C}$ does not belong to a  compact subgroup of $Sp$ and the
`compact' symplectic subspace  
$${\cal J}_{co}^{\bot, \C} = \oplus_{\rho: |\rho|= 1} {\cal
J}_{\gamma}^{\bot, \C}(\rho)\leqno(1.1a.3.co) $$ 
where $P_{\gamma}^{\C}$ does belong to a compact subgroup of $Sp.$ This
and the following decompositions are described in more detail in Klingenberg [Kl],
but also somewhat differently since in [Kl] the symplectic form is extended  to the
complexification as a sesquilinear form  rather than as a complex bilinear 
form.

The non-compact subspace has the further symplectic orthogonal decompositon
$${\cal J}_{nc}^{\bot, \C} = \bigoplus_{\rho: |\rho| < 1} [{\cal J}_{\gamma}^{\bot,
\C} (\rho)
\oplus {\cal J}_{\gamma}^{\bot, \C} (\rho^{-1})]\leqno(1.1a.4)$$
into symplectic complex  2-planes.  We may rewrite this decomposition in the form
$${\cal J}_{nc}^{\bot, \C} =
 {\cal J}_{s}^{\bot, \C} \oplus {\cal J}_{u}^{\bot, \C}\leqno(1.1a.5)$$
where
$${\cal J}_{s}^{\bot, \C} = \bigoplus_{\rho: |\rho| < 1} {\cal J}_{\gamma}^{\bot, \C}
(\rho),
\;\;\;\;\;\;\;
{\cal J}_{u}^{\bot, \C} = \bigoplus_{\rho: |\rho| > 1} {\cal J}_{\gamma}^{\bot, \C}$$
are the symplectically dual stable, resp. unstable Lagrangean subspaces. 
 
The compact subspace has the further symplectic decompositon
$${\cal J}_{co}^{\bot, \C} = 
\bigoplus_{\rho: \rho = e^{i\alpha}, \alpha \in (0, \pi)}
{\cal J}_{\gamma}^{\bot, \C}(\rho) \oplus {\cal J}_{\gamma}^{\bot, \C}(\bar{\rho}).
\leqno(1.1a.6)$$
Any choice of one $\rho$ from a pair $\{\rho, \bar{\rho}\}$ determines a splitting
of ${\cal J}_{co}^{\bot, \C}$ into a pair of dual Lagrangean subspaces.

On the level of real symplectic spaces, we have the closely related
$P_{\gamma}$-invariant symplectic decomposition
$$ {\cal J}_{\gamma}^{\bot} = {\cal J}_{s }^{r} \oplus
 {\cal J}_{ u}^{r} \oplus  {\cal J}_{
ce}^{\bot, 2p}\leqno(1.1a.7)$$ into the stable, unstable and center stable real
subspaces of dimensions
$r,r,2p$ respectively. By definition,
$$ {\cal J}_{s }^{r} = \bigoplus_{\rho \in \R, |\rho|<1}
 {\cal
J}^{\bot,\R}_{\gamma} (\rho) \oplus \bigoplus_{\rho \in \C - \R, |\rho|<1} 
 {\cal J}^{\bot,\R}_{\gamma}(\rho)\leqno(1.1a.8)$$
where 
$$P_{\gamma}|_{{\cal J}^{\bot,\R}_{\gamma}(\rho)} = \rho \;\;\;\;\;\;\;\;(\rho \in
\R)$$ 
 respectively
$$P_{\gamma}|_{{\cal J}^{\bot,\R}_{\gamma}(\rho)} = e^{-\mu}
\left( \begin{array}{ll} cos \nu & sin \nu \\ -sin \nu & cos \nu 
\end{array} \right)\;\;\;\;\;\;\;\;\;\;\;
\rho = e^{-\mu + i \nu}, \; \;\mu, \nu > 0. $$ In the latter case, 
${\cal J}^{\bot,\C - \R}_{\gamma}(\rho)$ is  the 
real symplectic 2-plane
whose complexification equals
$ {\cal J}_{\gamma}^{\bot, \C}(\rho) \oplus {\cal J}_{\gamma}^{\bot, \C}(\rho^{-1})$.
Similarly for the case of the unstable subspace. In the center stable case, 
$$ {\cal J}_{\gamma, ce}^{\bot, 2p} = \bigoplus_{j=1}^p {\cal J}^{\bot, \R}_{\gamma}
(\alpha_j)\leqno(1.1a.9)$$ where
$$P_{\gamma}|_{{\cal J}^{\bot,\R}_{\gamma}(\alpha)} = 
\left( \begin{array}{ll} cos \alpha & sin \alpha \\ -sin \alpha & cos \alpha 
\end{array} \right)$$
and with
${\cal J}^{\bot, \R}_{\gamma} (\alpha)$ the symplectic two plane whose
complexification equals
${\cal J}_{\gamma}^{\bot, \C}(\rho) \oplus {\cal J}_{\gamma}^{\bot, \C}(\bar{\rho})$
with $\rho = e^{i\alpha}.$

We will put $r = q + 2c$ where $q = \# \{\rho \in \R, |\rho|<1\}$ and where
$2c = \# \{\rho \in \C - \R, |\rho|<1\}$ and say that $\gamma$ has type $(p,q, c)$
if it has $p$ pairs of conjugate eigenvalues of modulus one
$\{ e^{i\alpha}, e^{-i\alpha}\},$
 $q$ pairs of real inverse eigenvalues $\{e^{\lambda}, e^{-\lambda}\}$, and
$c$ quadruples of non-real complex eigenvalues $\{e^{\pm \mu +  \pm i \nu}\}.$
Here, we have assumed the real eigenvalues are positive for brevity of notation.
  Finally we may write:
$${\cal J}^{\bot}_{\gamma} = {\cal J}^{e}_{\gamma} \oplus {\cal J}^{h}_{\gamma}
\oplus {\cal J}^{ch}_{\gamma}\leqno(1.1a.10)$$
where ${\cal J}^{e}_{\gamma} = {\cal J}_{\gamma, ce}^{\bot, 2p}$ is the {\it
elliptic} (or real center stable) subspace, where 
$${\cal J}^{h}_{\gamma} = \bigoplus_{\rho \in \R, |\rho|<1}
 {\cal J}^{\bot,\R}_{\gamma} (\rho) \oplus {\cal J}^{\bot,\R}_{\gamma}
(\rho^{-1})\leqno(1.1a.11)$$ is the {\it real hyperbolic} subspace and where 
$${\cal J}^{ch}_{\gamma} =  \bigoplus_{\rho \in \C - \R, |\rho|<1} 
 {\cal J}^{\bot,\R}_{\gamma}(\rho) \oplus {\cal
J}^{\bot,\R}_{\gamma}(\rho^{-1})\leqno(1.1a12)$$ is the {\it complex hyperbolic} (or
loxodromic) subspace.
\medskip

\noindent{\bf \S1.1b: Jacobi eigenvectors and Wronskian matrix}
\medskip

As mentioned in the introduction, the intertwining operator to the quantum normal
form will involve a certain Wronskian matrix of the Jacobi equation. Roughly
speaking, it is the real symplectic matrix whose entries are given by the real and
imaginary parts of the Jacobi eigenvectors and their time derivatives relative to
a normal frame.

More precisely, we fix a symplectic orthonormal basis of Jacobi eigenvectors as
follows:
\medskip

\begin{tabular}{l|l|l} complex subspace & eigenvectors & normalization  \\ \hline
elliptic plane & $P_{\gamma} Y_j^e = e^{i \alpha_j}Y_j^e, \;P_{\gamma} \bar{Y^e}_j =
e^{-i
\alpha_j}\bar{Y^e}_j $ & $\omega(Y_j^e, \bar{Y}^e_j) = 1.$\\ \hline
real hyp. plane & $P_{\gamma} Y^+_j = e^{\lambda_j}Y^{+}_j, \;P_{\gamma} Y^{-}_j =
e^{- \lambda_j}Y^{-}_j $ & $\omega(Y_j^+, Y^{-}_j) = 1.$\\ \hline
cx.hyp.4-plane & $P_{\gamma} Y^{++}_j = e^{\mu + i \nu }Y^{++}_j,\;
P_{\gamma} Y^{--}_j = e^{- \mu - i \nu}Y^{--}_j,$ & $\omega(Y^{++}, Y^{--}) = 1$ \\
\hline & $P_{\gamma} Y^{+-}_j = e^{ \mu - i \nu }Y^{+-}_j,\;
P_{\gamma} Y^{-+}_j = e^{- \mu + i \nu}Y^{-+}_j,$ & $\omega(Y^{-+}, Y^{+-}) = 1$\\
\end{tabular}

The normalization makes sense since complex 2-planes spanned by eigenvectors
corresponding to inverse eigenvalues are symplectic.

We now fix   a parallel
normal frame $e(s):= (e_1(s),...,e_n(s))$ along $\gamma|_{[O,L)}$ and denote by
$\langle Y, e_j \rangle$ the Riemannian inner product of a vector  $Y$ along
$\gamma$ with the the jth normal vector. Corresponding to the splitting
of ${\cal J}^{\bot}_{\gamma}$ into its elliptic, hyperbolic and and complex
hyperbolic (real) subspaces, we then get a real symplectic 2n x 2n Wronskian matrix  
$${\cal W} = [{\cal W}^e |{\cal W}^h | {\cal W}^{ch}]\leqno(1.1b.1)$$
formed by the 2n x 2p elliptic Wronskian matrix
$${\cal W}^e(s) := \left ( \begin{array}{ll}  Re\langle Y^e_i, e_j\rangle  &
Im \langle Y_i^e ,e_j \rangle  \\  Re \langle \dot{Y_i}^e,e_j \rangle  &  Im \langle
\dot{Y^e_i}, e_j \rangle  \end{array} \right )_{i=1,\dots,p;
j=1,\dots,n}\leqno(1.1b.2e)$$ the 2n x 2q real  hyperbolic Wronskian matrix
$${\cal W}^h(s) := \left ( \begin{array}{ll} \langle Y_i^+, e_j \rangle 
 &  \langle Y_i^- ,e_j \rangle  \\\langle \dot{Y_i}^+,e_j \rangle & \langle
\dot{Y_i}^-,e_j \rangle \end{array} \right
)_{i =1,\dots,q;j=1,\dots,n}\leqno(1.1b.2h)$$ and finally the complex hyperbolic 2n x
4c Wronskian matrix
$${\cal W}^{ch}(s) := \left ( \begin{array}{llll}  Re \langle Y_i^{++},e_j
\rangle & Im \langle Y_i^{++}, e_j \rangle & Re \langle Y_i^{--},e_j
\rangle & Im \langle Y_i^{--}, e_j \rangle\\
  Re \langle \dot{Y_i}^{++},e_j
\rangle & Im \langle \dot{Y_i}^{++}, e_j \rangle & Re \langle \dot{Y_i}^{--},e_j
\rangle & Im \langle \dot{Y_i}^{--}, e_j \rangle \end{array}
\right)_{i=1,\dots,c; j=1,\dots,n}.\leqno(1.1b.2ch)$$
\medskip

To see that ${\cal W}$ is indeed symplectic, and 
to better understand its properties, we reconsider
the Jacobi equation and Poincare map from the Riemannian viewpoint.  Thus, we let
$\nabla$ denote the Riemannian connection, and recall that it determines a horizontal
subbundle of $T(S^*M)$ complementary to the vertical subbundle of the projection
$\pi: S^*M \rightarrow M.$  Together with the symplectic structure, we get a
splitting $$T(S^*M) = \bar{H} \oplus  \bar{V} \oplus \bar{T}$$
where $\bar{T}$ is the real span of $\dot{\gamma}$, and $\bar{H} \oplus  \bar{V}$
is the horizontal plus vertical decomposition of the kernel of the contact form
$\alpha = \xi \cdot dx$ (or equivalently, of the symplectic orthogonal of $T$ and
the cone axis).  The subspaces $\bar{H},  \bar{V}$ are symplectically paired
Lagrangean subspaces of $T(T^*M)$. Given a vector $X \in N_{\gamma(t)}$, we denote by
$X^h$ the horizontal lift of $X$ to $\bar{H}_{\gamma(t)}$ and by $X^v$ the vertical
lift to $\bar{V}_{\gamma(t)}$.  The correspondence
$$Y(t) \rightarrow (Y(t)^h, \dot{Y}(t)^v)$$
then defines an isomorphism between the spaces of Jacobi fields and geodesic flow
invariant vector fields along $(\gamma(t), \dot{\gamma}(t)$ (cf [Kl, Lemma 3.1.6]). 
That is, 
$$dG_{(\gamma(0), \dot{\gamma}(0))}^s: (Y(0)^h, \dot{Y}(0)^v) \rightarrow
(Y(s)^h, \dot{Y}(s)^v)$$
where $Y(s)$ is the Jacobi field with the given initial conditions.  Moreover, since
$G^t$ is a Hamiltonian flow $dG^s$
is a linear symplectic mapping from $(\bar{H} \oplus  \bar{V})_{\gamma(0),
\dot{\gamma}(0)}$ to $(\bar{H} \oplus  \bar{V})_{\gamma(s), \dot{\gamma}(s)}.$

Relative to the basis $\{e_j(s)\}$,
 the Jacobi equation is equivalent to the linear system
$\frac{D}{ds}(Y,P) = J H (Y, P)$, where: (i) $P=\frac{DY}{ds},$ (ii) $J$ is the
standard  complex structure on $R^{2n},$ and where (iii)
$$H =\left (\begin{array}{ll} K & 0 \\ 0 & I \end{array} \right ) $$
with $K$  the curvature matrix and with $I$  the identity matrix. Moreover
the basis $\{e_j(s)\}$ induces a  moving symplectic frame  $\{e^h_j(s),
e^v_j(s)\}$ of $(\bar{H} \oplus\bar{V})_{\gamma(s), \dot{\gamma}(s)}.$  The evolution
operator for the linear system is then just $dG^s$ expressed as a matrix relative to
the moving symplectic frame. 

Now consider the above
 basis  $\{Re Y_j^e, Im Y_j^e, Y_j^+, Y_j^-, Re Y_j^{++},
Im Y_j^{++}, Re Y_j^{--}, Im Y_j^{--}\}$ of ${\cal J}^{\bot}_{\gamma}$ formed by 
the eigenvectors of $P_{\gamma}$. By construction, it is a symplectic basis  relative
to the Wronskian form $\omega$. Hence the pairs $(Y^h(s),P^v(s))$ consisting of these
eigenvectors and their time derivatives form a moving symplectic basis of 
$\bar{H} \oplus \bar{V}.$  Expressed in terms of the  frame
$\{e_j^h(s), e_j^v(s)\}$ we then get a symplectic basis of $\R^{2n}$ relative
to the standard symplectic structure. The Wronskian matrix ${\cal W}(s)$ is just
the matrix whose columns are formed by these basis elements, and it is therefore
symplectic for each $s$.

Consider now the monodromy aspect of ${\cal W}(s)$, i.e. its transformation law
under time translation $s \rightarrow s+L$ thru one period. Let $Y_i$ denote one
of the complex eigenvectors of $P_{\gamma}$. Then the  
 matrix element $\langle Y_i(s),e_j(s)
\rangle$ (or with $\dot{Y}_i$ in place of $Y_i$) satisfies
$$\langle Y_i(s+L),e_j(s+L) \rangle = \rho_i \sum_{k=1}^n t_{jk} \langle Y_i(s), 
e_k(s) \rangle\leqno(1.1b.3)$$
where $T:= (t_{jk})$ is the holonomy matrix,
$$e_j(s+L) = \sum_{k=1}^n t_{jk} e_k(s).$$ 
It follows that
$${\cal W}(s+ L) = {\cal W}(s) T^* P_{\gamma}. \leqno(1.1b.4)$$
\medskip

\noindent{\bf \S1.1c: Symplectic equivalence of quadratic Hamiltonians}
\medskip

Let $(\R^{2n}, \omega)$ be the standard symplectic vector space, endowed with
linear coordinates $x= (q_1,\dots,q_n, p_1,\dots, p_n)$ such that
$\omega = \sum dq_i \wedge dp_i.$  A quadratic (real) Hamiltonian is by definition a
quadratic form
$$H(q,p) = \half \langle Ax,x\rangle = \half \omega( JA x, x)$$
where $\langle \cdot \rangle$ is the Euclidean scalar product, where $A$ is
a 2n x 2n real symmetric matrix and where $J= \left( \begin{array}{ll} 0 & -I \\
I & 0 \end{array} \right).$  Then $JA \in sp(\R^{2n}, \omega)$, and hence its
spectrum decomposes into purely imaginary
pairs $(i\alpha, -i\alpha)$, into
 real pairs $(\lambda, -\lambda)$, and into complex quadruples $(\pm \mu \pm
i\nu).$ We will assume, as above, that the eigenvalues are simple and not equal to
$\pm 1.$ Then $H(q,p)$ decomposes into sums of terms of the following types
of quadratic Hamiltonians, or classical `actions':
the elliptic type
$$I^e (q_1,p_1) := \half \alpha (q_1^2 + p_1^2), \leqno(1.1c.1e)$$
 the real hyperbolic type
$$I^{h}(q_1,p_1):= \half \lambda q_1 p_1, \leqno(1.1c.1h)$$
and the complex hyperbolic (or loxodromic) type
$$I^{ch}(q_1,p_1,q_2, p_2) =  \half \mu (q_1p_1 + q_2p_2) +
\half \nu ( q_1 p_2 - q_2 p_1 ). \leqno(1.1c.1ch)$$
Note that
$$q_1p_1 + q_2p_2 = Re (q_1 + i q_2) ( p_1 - i p_2),\;\;\;\;\;\;\;
 q_2 p_1 - q_1 p_2  = Im (q_1 + i q_2) ( p_1 - i p_2)$$ and
$$\{q_1p_1 + q_2p_2, q_1 p_2 - q_2 p_1\} = 0$$ 
To unify these expressions, we observe that they all have the form $\half Re\; s a^*
a$  where $s \in \C$ and where 
$a^*, a$ denote symplectically dual complex linear coordinates.  Indeed, in the
 elliptic case, $a = q_1 + i p_1, a^* = q_1 - i p_1, s= \alpha \in \R;$ in the real
hyperbolic case, $a = q_1, a^* = p_1, 
s= \lambda \in \R;$ and in the loxodromic case,
$s = \mu + i\nu, a = (q_1 + i q_2), a^* = (p_1 - i p_2).$

\medskip

\noindent{\bf \S1.2: Microlocal preliminaries}
\medskip

As mentioned above,  
 the  wave invariants only involve the metric and Laplacian 
$\Delta$ in a tubular neighborhood of $\gamma$.
 In fact, as discussed in [G.1][Z.1] they only involve 
the microlocalization
 of $\Delta$ to the conic neighborhood
$$|y| \leq \epsilon, \;\;\;\;\;\;\;\; |\eta| < \epsilon \sigma \leqno(1.2.1)$$
 of $T^*S_L^1-0$ in
$T^*(S^1 \times \R^n)$.
Here,$(s,\sigma,y,\eta)$ denote the symplectic Fermi coordinates
and $\psi$ denotes a smooth homogeneous cut-off function on
$T^*(S_L^1 \times \R^n)-0$ which equals 1 in some conic neighborhood $V$ of
$T^*S_L^1-0$ and vanishes identifically off of some slightly larger conic
neighborhood.  

As in the case of elliptic closed geodesics, 
 to put $\Delta$ into a microlocal
(quantum Birkhoff) normal form around $\gamma \sim S^1$ is first of all to
conjugate it to a distinguished maximal abelian subalgebra ${\cal A}_{\gamma}$
of the algebra $\Psi^*(S^1_L \times \R^n)$ of pseudo-differential operators on
the model space $S_L^1 \times \R^n$. This distinguished subalgebra will depend
on the type of the geodesic $\gamma$. Roughly speaking, it will consist of
the tangential operator $D_s:= \frac{1}{i} \frac{\partial}{\partial s}$
together with an appropriate set of quantized quadratic normal forms or
`action operators' in the
transversal directions. 
\medskip

\noindent{\bf \S1.2.1: The model algebras}
\medskip

To specify this `appropriate set' of action operators, we begin by recalling that 
 the Schrodinger representation of 
the (complexified) Heisenberg algebra ${\bf h}_n\otimes \C$ on the transverse
space $L^2(\R^n)$, is generated by
the self-adjoint operators
$Y_j =$ ``multiplication by $y_j$" and 
by $D_j=\frac{\partial}{i\partial y_j}$.
Equivalently it is generated by the creation/annihilation operators
$Y_j + i D_j, Y_j - i D_j.$
For our purposes, however, it will be more natural to choose a different set of
generators depending on the $(q,p,c)$ type of the closed geodesic. 

Corresponding to
the $2p$ dimensional elliptic symplectic subspace we will use
as generators  the above (elliptic) annihilation/creation operators
$$Z_j := Y_j + iD_{y_j}  \;\;\;\;\;\;\;\; Z_j^{\dagger}
 = Y_j - i D_{y_j}\;\;\;\;\;(j=1,\dots,
p)\leqno(1.2.1.1e)$$ which satisfy the commutation relations
$$[Z_j, Z_k] = [Z_j^{\dagger}, Z_k^{\dagger}]
 = 0 \;\;\;\;\;\; [Z_j, Z_k^{\dagger}] =
2\delta_{ij} I.$$ 

We would like to use the real (resp. complex) hyperbolic analogues in the
hyperbolic subspaces.  To determine the analogues
we note that $Z_j$, resp. $Z_j^{\dagger}$ 
are the Weyl quantizations
of the  symplectically dual (modulo a factor of 2) complex linear coordinates
$z_j:= y_j + i \eta_j,$ resp. $z_j^{\dagger}:= y_j - i \eta_j$.  
We use the `dagger' notation rather than the adjoint notation $Z_j^*$ to
emphasize that the dual operators are symplectically dual; they are also
adjoints of each other, but this property will not extend to the hyperbolic
cases.
Indeed, corresponding to the $2q$ dimensional real hyperbolic subspace, the
natural generators are the hyperbolic annihilation/creation operators
$$ Y_j ,\;\;\;\;\;\;\;\;\;\; D_{y_j} \;\;\;\;\;(j=p + 1,\dots,q)\leqno(1.2.1.1h)$$
which of course are also symplectically dual.  And corresponding to the $4c$
dimensional complex hyperbolic subspace, we the natural generators are the
complex hyperbolic annihilation/creation operators
$$W_j := Y_j + i Y_{c + j},\;\;\;\;\;W_j^{\dagger}:= D_{y_j} - i
D_{y_{c+j}},\;\;\;\;\;\bar{W_j}:= Y_j - i Y_{c+j},\;\;\;\;\;\;\bar{W_j}^{\dagger}:=
D_{y_j} + i D_{y_{c+j}} \;\;\;\;\;(j= p + q + 1,\dots,c).\leqno(1.2.1ch)$$
We note that they satisfy the commutation relations:
$$[W_j, W_j^{\dagger}] = -2,\;\;\;\;\;\;\;\; [\bar{W_j},\bar{W_j}^{\dagger}] = -2$$
with all other brackets zero.  We will not bother to renormalize the operators to
be precisely dual.

The enveloping algebra of the Heisenberg algebra is then generated by all the
above annihilation/creation operators, 
$${\cal E} := <Z_1,\dots Z_p, Z^{\dagger}_1,\dots, Z_p^{\dagger},Y_1,...,Y_q,
D_{y_1},...,D_{y_q}, W_1,\bar{W}_1,\dots, W_c,\bar{W}_c,
W_1^{\dagger},\bar{W}_1^{\dagger},
\dots,W_c^{\dagger},\bar{W}_c^{\dagger}>\leqno(1.2.1.2)$$ 
 and is of course 
the algebra of partial differential operators on $\R^n$ 
with polynomial coefficients. We will denote
by 
${\cal E}^n$ the subspace of polynomials of degree n 
in the generators. The microlocalization of this algebra 
is  the  isotropic Weyl algebra ${\cal W}^*$ of 
pseudo-differential operators on $\R^n$, in which the
generators are assigned the order $\frac{1}{2}$, so that
$${\cal E}^n \subset {\cal W}^{n/2}$$
$$[{\cal E}^m, {\cal E}^n] \subset {\cal E}^{m + n - 2}. $$
The symplectic algebra ${\bf sp}(n,\C)$ is then represented in ${\cal E}^2$ by 
homogeneous quadratic polynomials in the generators, which have degree 1. In
particular it contains the following elliptic, resp. hyperbolic, resp. complex
hyperbolic ( loxodromic)
`action' operators: 
$$\hat{I}^e_{j}:= Z^{\dagger}_j Z_j, \;\;\;\;\;\;\;\hat{I}^h_j = \half (Y_j D_{y_j} +
D_{y_j}Y_j),\;\;\;\;\;\;\;\hat{I}^{ch,Re}_j = \half Re(
W_j^{\dagger}W_j + (W_j^{\dagger}W_j)^*),\;\;\;\;\;\hat{I}^{ch,Im}_j = Im
W_j^{\dagger}W_j.\leqno(1.2.1.3)$$
The complex hyperbolic action operators can also be written in the form
$$\hat{I}^{ch,Re}_j = \half(Y_j D_{y_j} + D_{y_j}Y_j + Y_{j+c}D_{y_{j+c}} +
D_{y_{j+c}}Y_{j+c}),\;\;\;\;\;\;\hat{I}^{ch,Im}_j = (Y_j D_{y_{j+c}} -
Y_{j+c}D_{y_j})\leqno(1.2.1.4)$$ where the coordinates are indexed so that the
$dy_j\wedge dy_{j+c}
\wedge d\eta_j \wedge d\eta_{j + c}$-planes are  the $P_{\gamma}$-invariant complex
hyperbolic 4-planes.  It is then natural to introduce polar coordinates $r_j, \phi_j$
on the $(y_j, y_{j+c})$-plane so that the loxodromic actions operators simplify to
$$I^{ch,Re}_j = \half (r_j D_{r_j} + D_{r_j} r_j),\;\;\;\;\;\;\;I^{ch,Im} =
D_{\theta}.\leqno(1.2.1.5)$$

  We now introduce the  distinguished $(p,q,c)$
 maximal (transverse) abelian subalgebra of ${\cal W}$,  given by 
$${\cal I}_{p,q,c}:= <I^e_1,...,I^e_p, I_1^h,\dots, I_q^h,
I_1^{ch,Re},\dots,I_{c}^{ch,Re}, I_1^{ch,Im},\dots, I_c^{ch,Im}>.\leqno(1.2.1.6) $$
Together with the tangential operator we get the
{\it (p,q,c)- maximal abelian subalgebra } given by
$${\cal A}_{p,q,c}:= <D_s,I^e_1,...,I^e_p, I_1^h,\dots, I_q^h,
I_1^{ch,Re},\dots,I_{c}^{ch,Re}, I_1^{ch,Im},\dots, I_c^{ch,Im}>.\leqno(1.2.1.7)$$
\medskip

\noindent{\bf \S 1.2.2: Model eigefunctions}

 An orthonormal basis of $L^2(S^1_L\times \R^n)$ of joint eigenfunctions of
 ${\cal A}_{p,q,c}$ is given as follows: corresponding to the (p,q,c)- type
of $P_{\gamma}$ we can write $$L^2(S^1_L\times \R^n) = L^2(S^1_L) \otimes L^2(\R^p)
\otimes L^2(\R^q) \otimes L^2(\R^{2c})$$ and construct the eigenfunctions as
(tensor) products of the eigenfunctions on the factors.  In the elliptic factors,
the eigenfunctions are the normalized Hermite functions $\gamma_q$ (cf.[F], or [Z.1]
for a context similar to the one here).  In the real hyperbolic factors,the
action operators are the generators of the unitary dilations on $L^2(\R)$ given by
$$U(\theta) f(x) = \theta^{\half} f(\theta x), \;\;\;\;\;\;\;(\theta \in \R^+).$$
Their generalized eigenfunctions are the temperate distributions
$$x_{+}^{-\half + ia},\;\;\;\;\;\;\;x_{-}^{-\half + ia},\;\;\;\;\;(a \in \R)$$
and any $f \in L^2(\R,dx)$ has the eigenfunction expansion
$$f(x) = \int_{\R} \hat{f}_{+}(a) x_{+}^{-\half + ia} da + \int_{\R} 
\hat{f}_{-}(a) x_{-}^{-\half + ia} da$$
with $\hat{f}_{\pm a}:= \langle f, x_{\pm}^{-\half + ia} \rangle.$
In the complex hyperbolic (i.e. loxodromic) factors, the actions operators are given
by the unitary dilations in polar coordinates
$$U(\rho) f(r, \theta) = \rho f(\rho r, \theta), \;\;\;\;\;\;\;(\rho \in \R^+)$$
together with rotations.  The joint eigenfunctions are  the temperate
distributions on $\R^2$ given by
$$r^{it -1} e^{in \theta}, \;\;\;\;\;\;\;(t\in \R)$$ 
and in a notation similar to that of the real hyperbolic case a function 
$f\in L^2(\R^2, r dr d\theta)$ may be expressed in the form
$$f(r,\theta) =  \sum_{ n \in \Z} \int_{\R} \hat{f}(t,\theta) r^{it -1}e^{in \theta}
dt.$$
For future reference we summarize the situation in the following table.

\medskip

\begin{tabular}{l|l|l} Factor & Action & Eigenfunction \\ \hline
$L^2(S^1_L)$ & $D_s$ & $e^{is \frac{2\pi k}{L}}$ \\ \hline 
$L^2(\R^p)$ & $I_j^e$ &    Hermite functions  $\gamma_{ m}, m \in \Nb^n;$\\
$\cdot $ & $\cdot$ &  $\gamma_o(y)=\gamma_{iI}(y):= e^{-\half |y|^2}$\\
$\cdot $ & $\cdot$ &
$ \gamma_{ m}:= C_m a_1^{\dagger m_1}...a_n^{\dagger m_n} \gamma_o $ \\ \hline
$L^2(\R^q)$ & $I_j^h$ & $\Pi_{j=1}^r y_{j \pm}^{i a_j - \half}, a \in \R^r$\\ \hline
$L^2(\R^{2c})$ & $I_j^{ch Re}, I_j^{ch, Im}$ & $\Pi_{j=1}^{c} r_j^{ i t_j - 1}
e^{i n_j \theta_j}, t \in \R^{c}$ \\ \hline 
\end{tabular}
  \medskip

As in the introduction, we put:
$$H_{\alpha,\lambda,(\mu,\nu)} :=  [\sum_{j=1}^p \alpha_j \hat{I}_j^e + 
 \sum_{j=1}^q \lambda_j
\hat{I}_j^h + \sum_{j=1}^c \mu_j \hat{I}_j^{ch,Re} + \nu_j
\hat{I}_j^{ch,Im}]\leqno(1.2.2.1)$$
$${\cal R}:= \frac{1}{L}(L D_s +  H_{\alpha,\lambda, (\mu,\nu)})$$
and note that
$${\cal R} e^{is \frac{2\pi k}{L}}\gamma_{ m}(x)[\Pi_{j=1}^r y_{j \pm}^{i
a_j - \half}][\Pi_{j=1}^{c} r_j^{ i t_j-1}e^{i n_j \theta_j}] = r_{k m n a t} 
e^{is \frac{2\pi k}{L}}\gamma_{ m}(x)[\Pi_{j=1}^r y_{j \pm}^{i
a_j - \half}][\Pi_{j=1}^{c} r_j^{ i t_j-1}e^{i n_j \theta_j}]\leqno(1.2.2.2a)$$
with
$$r_{k m n a t} = \frac{1}{L}(2 \pi k + [\sum_{j=1}^p \alpha_j (m_j + \half) + 
\sum_{j=1}^r \lambda_j a_j  + \sum_{j=1}^c \mu_j t_j + \nu_j n_j]).\leqno(1.2.2.2b)$$
\medskip

\section{ The semi-classically scaled Laplacian}

  The significance of the maximal abelian algebra ${\cal A}_{pqc}$ will appear as soon
as we semi-classically `rescale' the Laplacian and conjugate the principal part,
the `linearized  $\sqrt{\Delta}$', to its normal form.

Let us briefly recollect this rescaling, which proceeds exactly as in the purely
elliptic case [Z.1,\S 2].  We first prepare the Laplacian by putting it in Fermi
normal coordinates $(s,y)$.  It is then self-adjoint relative to the volume
density  $J(s,u)|ds||dy|$ in these coordinates. To simplify, we then
 conjugate it to the unitarily equivalent  (1/2-density-)   Laplacian
$$\Delta_{1/2} := J^{1/2} \Delta J^{-1/2},$$ 
which is self-adjoint with respect to the Lesbesgue density $|ds dy|$.  

We thus have:
$$-\Delta_{1/2} = J^{-1/2}\partial_s g^{oo}J \partial_s J^{-1/2}
+\sum_{ij =1}^{n} J^{-1/2}\partial_{y_i} g^{ij} J \partial_{ y_j}
J^{-1/2}\leqno(2.1.1)$$
$$\equiv g^{oo}\partial_s^2 + \Gamma^o \partial_s + 
 \sum_{ij=1}^n g^{ij} \partial_{u_i}\partial_{y_j} + \sum_{i=1}^{n} \Gamma^{i}
\partial_{y_i} + \sigma_o.$$

Semi-classical rescaling then involves two conjugations: First, by 
$M_h =$ multiplication by $e^{\frac{is}{Lh}}$,
$$-M_h^* \Delta M_h = -(hL)^{-2}g^{oo} + 
2i(hL)^{-1}g^{oo} \partial_s + i(hL)^{-1}\Gamma^o +
\Delta $$
and then by  the semi-classical dilation $T_h f(s,y) =
f(s, h^{-\half}y)$.   
The complete conjugation $-T_h^*M_h^* \Delta M_h T_h$  results in the {\it
semi-classically scaled Laplacian}
$$-\Delta_{h} =  -(hL)^{-2} g^{oo}_{[h]}
 + 2i(hL)^{-1}g^{oo}_{[h]}\partial_s + i(hL)^{-1}\Gamma^o_{[h]}+
 h^{-1}( \sum_{ij=1}^n g^{ij}_{[h]}\partial_{y_i}\partial_{y_j})
 + h^{-\half}(\sum_{i=1}^{n} \Gamma^{i}_{[h]}
\partial_{y_i}) + (\sigma)_{[h]},\leqno(2.1,2)$$
 the subscript $[h]$ indicating to dilate
 the coefficients of the operator in the form,
$f_h(s, y):=f(s, h^{\half} y).$  

Expanding the coefficients in Taylor series at $h=0$, we obtain 
the semi-classical  expansion
$$\Delta_h \sim \sum_{m=0}^{\infty} h^{(-2 +m/2)}{\cal L}_{2-m/2}
 \leqno (2.1.3)$$
where ${\cal L}_2 = L^{-2},$ ${\cal L}_{3/2}=0$ and where
$${\cal  L}_1 = 2 L^{-1}[i  \frac{\partial}{\partial s} + \half \{\sum_{j=1}^{n}
\partial_{y_j}^2 -
\sum_{ij=1}^{n} K_{ij}(s) y_i y_j\}] \leqno (2.1.4).$$
We will denote the bracketed operator, the `linearized $\sqrt{\Delta}$' by ${\cal
L}$.  It is of order 1 in the sense of pseudodifferential operators (using the Weyl
filtration in the transverse variables) and as will be seen below is the principal
term in the semi-classical expansion of the square root of the rescaled Laplacian.
\bigskip

\noindent{\bf (2.1.A) Appendix on metric scaling}
\bigskip

In addition to semi-classical scaling, we have also just introduced an independent
scaling, {\it metric scaling}, which has to do with the behaviour of objects under
dilations $g \rightarrow \epsilon^2 g$ of the metric.  As discussed in detail in 
[Z.1], the wave invariants have well-defined weights under metric rescaling and
in analysing them it is very convenient to rescale all objects to be weightless.
For instance,  as discussed in
[Z.1, \S 1.4],  an $\omega$- symplectic basis of Jacobi fields has weight $\half$ and
its time derivative has weight $-\half.$  To render it weightless  a Jacobi
eigenfield
$Y$ should be replaced by $L^{- \half} Y$, $\dot{Y}$ by $L^{\half} \dot{Y}$ etc. The
resulting weightless Wronskian matrix is then denoted by ${\cal W}_L$. It is
essentially the weightless matrix denoted ${\cal A}_L$ in [Z.1].

To render the coordinates $(y, \eta)$ weightless under metric rescaling,  
we also change variables to $x = L^{-1} y$ and rewrite
$\Delta_h$ and the ${\cal L}_{2-\frac{n}{2}}$'s in terms of the $x$-variables.
For instance,  ${\cal L}$ then takes the form: 
$${\cal L} = i  \frac{\partial}{\partial s} +
 \half [\sum_{j=1}^{n}L^{-1} \partial_{x_j}^2 -
\sum_{ij=1}^{n} L K_{ij}(s) x_i x_j].$$
The symplectic coordinates on the symplectic normal space $T^*\R^n$ to $\R^+ \gamma$
will henceforth be denoted $(x, \xi).$

For the sake of brevity we will not  draw much attention to metric scaling in the
various steps to come in the normal form algorithm. In all cases, the role of metric
scaling is identical to that in the elliptic case [Z.1].   
\medskip

\noindent{\bf \S2.2: Conjugating $\Delta_h$ to the model}
\medskip

We now conjugate the semi-classically scaled Laplacian $\Delta_h$ from
$L^2(N_{\gamma})$ to the model $L^2(S^1_L \times \R^n)$ by means of the moving
metaplectic operator $\mu({\cal W}_L)$,
$$\mu({\cal W}_L) f(s,y) :=  \mu({\cal W}_L(s))f(s,y).$$ 
The motivation for this conjugation comes from:
\medskip

\noindent{\bf (2.2.1) Proposition}{\it  
$${\cal L} = \mu({\cal W}_L^*) D_s \mu({\cal W}_L^*)^{-1}$$
where ${\cal W}_L$ is the weightless Wronskian matrix and $\mu$ is the metaplectic
representation.}
\medskip

\noindent{\bf Proof}: 

First, let us ignore the scaling parameter $L$, i.e. let us put $L=1$.
 The right side is then equal to
$(D_s + \mu({\cal W}(s))^{*} D_s\mu( {\cal W}(s)).$ 
  To evaluate the second term, we recall that the columns of ${\cal W}$
 are Jacobi fields, and that Jacobi's equation is equivalent to the linear system
$\frac{D}{ds}(Y,P) = J H (Y, P)$  with $P=\frac{DY}{ds},$ and with
$$H =\left (\begin{array}{ll} K & 0 \\ 0 & I \end{array} \right ). $$
  Hence, the second term is $\frac{1}{i} d\mu (JH)$ with
$d\mu$ the derived metaplectic representation.   But $\frac{1}{i}d\mu(JH) =
1/2(\sum_{i=1}^{n} \partial_{y_i}^2 - \sum_{ij=1}^{n} K_{ij}(s)y_i y_j)$ [F]. 
Re-inserting $L$ to make all objects weightless, we get the formed claimed above. 
\qed
\medskip

 Thus, 
conjugation by ${\cal W}_L$ puts the principal term ${\cal L}$ of $\Delta_h$ into the
simple normal form $D_s$.  This suggests conjugating the  full rescaled Laplacian by
$\mu({\cal W}_L)$ to  the `twisted model'
semi-classical Laplacian
$${\cal D}_h =\mu({\cal W}_L^*)^{-1} \Delta_h \mu({\cal W}_L^*)\leqno(2.2.2)$$
which  has the asymptotic expansion
$${\cal D}_h \sim \sum_{m=o}^{\infty} h^{(-2 + \frac{m}{2})}
{\cal D}_{2 -\frac{m}{2}} \leqno(2.2.3)$$
  with  ${\cal D}_2 = I, {\cal D}_{\frac{3}{2}}=0,
{\cal D}_1=D_s$. Thus, ${\cal D}_h$ is a small perturbation of $D_s$,
and one may expect that  perturbation theory can be used to find a good
normal form for the whole of ${\cal D}_h$. 

Before doing so, we must
 consider which Hilbert space is the natural domain for ${\cal D}_h$. The point is
that the conjugation has non-trivial monodromy (\S1.1b) and hence the conjugate will
act on functions transforming correctly under the monodromy group.  

We can describe the Hilbert space in terms of quantum mapping cylinders [Z.1].
First, we consider the holonomy aspect, put
$$C^{\infty}_T(\R\times \R^n):= \{ f \in C^{\infty}(\R \times \R^n):
f(s+L,u) = \mu(T) f(s,u)\}\leqno (2.2.4)$$
and let ${\cal H}_T$ denote its closure with respect to the obvious inner product
over $[0,L)\times \R^n$.  Note that the metaplectic operator $\mu(T)$ is simply
$$\mu(T) f(u) = f(t^{-1}u)$$
and hence that  
$$C^{\infty}_T(\R \times \R^n)\sim C^{\infty} (N_{\gamma})$$
where the isomorphism is simply the pull-back by the exponential map defined by
the frame $e(s)$. In other words, expressed in terms of Fermi
coordinates relative to a normal frame, $L^2(N_{\gamma})$ becomes the quantum
mapping cylinder of $\mu(T^*).$
 Let us note however that $\Delta_h$ and
hence all the ${\cal L}_{2 -\frac{k}{2}}$'s are invariant under under $\mu(T)$, 
so that it will play an insignificant role for our purposes.

 On the other hand, the 
quantized linear Poincare map $\mu(P_{\gamma})$ will play an essential role. Hence
we introduce its quantized mapping cylinder
$${\cal H}_{\gamma} := \{ f \in L^2_{loc} (\R \times \R^n): \tau_L f =
\mu(P_{\gamma}) f \}\leqno(2.2.5)$$
and note that 
$$\mu({\cal W}_L ) : L^2(N_{\gamma}) \rightarrow {\cal H}_{\gamma}$$
is a unitary equivalence.
Hence, the natural domain for ${\cal D}_h$ is the quantum mapping cyliner of
$\mu(P_{\gamma})$.

In the calculation of traces, it is simpler to work in the original model
$L^2(S^1_L \times \R^n)$.  Hence we will also consider the conjugate of 
${\cal D}_h$ under a conjugation which untwists the mapping cylinder of
$\mu(P_{\gamma}).$   That is, we connect
$P_{\gamma}$ to the identity by a segment of the  one-parameter subgroup
$P_{\gamma}(s)$ thru $I$ and
$P_{\gamma}$, which exists by our non-degeneracy assumption on $P_{\gamma}$. Indeed,
after diagonalizing $P_{\gamma}$ and consulting the list of   symplectic
equivalence classes  of quadratic forms, we see that
$$P_{\gamma} = exp( \Xi_{H_{\alpha, \lambda, (\mu, \nu)}})$$
where $\Xi_f$ denotes the Hamilton vector field of $f$ and where $exp \circ \Xi$
denotes the exponential map from $sp(n,\R)\rightarrow Sp(n,
R)$, with $sp(n,\R)$ viewed as the Poisson
algebra of  quadratic functions on $\R^{2n}$.    
We then  have $$P_{\gamma}(s) = exp( s \Xi_{H_{\alpha, \lambda, (\mu, \nu)}})$$
and quantize this subgroup as
$$\mu(P_{\gamma}(s)) = e^{i s H_{\alpha, \lambda, (\mu, \nu)}} = e^{i s [\sum_{j=1}^p
\alpha_j I_j^e + 
\sum_{j=1}^q
\lambda_j I_j^h + \sum_{j=1}^c \mu_j I_j^{ch,Re} + \nu_j
I_j^{ch,Im}]}.\leqno(2.2.6)$$

Conjugation by $\mu(P_{\gamma})$ transforms ${\cal D}_h$
 into the model semi-classically scaled Laplacian
$${\cal R}_h := \mu(P_{\gamma} ) \mu({\cal W}_L) \Delta_h
\mu({\cal W}_L)^* \mu(P_{\gamma} )^* \sim
\sum_{m=o}^{\infty} h^{(-2 +
\frac{m}{2})} {\cal R}_{2 -\frac{m}{2}}\leqno(2.2.7) $$
with ${\cal R}_2 = I, {\cal R}_{\frac{3}{2}}=0,$ and with $ {\cal R}_1 := {\cal R}.$
All coefficients of terms in ${\cal R}$ are
periodic in $s$ and have weight -2 under metric rescaling.

\section{Semi-classical normal form}

We now wish to put ${\cal R}_{h}$
into semi-classical normal form, in the sense of [Z.1, Lemma 2.22]. This is
the key transitional step in putting $\Delta$ into microlocal normal form and
is the source of the connections to local geometric invariants.
 The method
is essentially the same as in the elliptic case, both in method and  in detail.
  We therefore present only
the first two steps in the proof and refer to [Z.1, loc.cit] for the inductive
argument.

As in the elliptic case, we state the result  in terms of the ${\cal
R}$-operators since  the trace will later be analysed in this  model. 
However, most of the proof will take place in the twisted model, where the
`linearized Laplacian' is $D_s$ and the equations
simplify most. In the following, the notation $|_o$ means to restrict to functions in
the kernel of ${\cal R}$, that is, to elements of ${\cal R}$-weight zero.  In
the twisted model, these are simply functions independent of $s$.  In the passage from
the semi-classical normal form to the microlocal (quantum Birkhoff) normal form,
the various operators will only be applied to such weightless elements.  This 
explains the rather complicated statement to follow; the result is only simple
and natural when restricted to elements of weight zero.
\medskip

\noindent(3.1)~~~~~~~{\bf Lemma (cf. [Z.1, Lemma 2.22]}~~~~~~~~~ {\it There exists an
$L$-dependent
$h$-pseudodifferential operator
$W_h=W_h(s, x, D_x)$ on $L^2(S^1_L \times \R^n)$ such that,
for each $ s\in S^1_L$, $$W_h(s,x,D_x): L^2(\R^n)\rightarrow L^2(\R^n) $$
is unitary, and such that
$$W_h^* {\cal R}_{h} W_h \sim - h^{-2}L^{-2}
 + 2 h^{-1}L^{-1}{\cal R} +
\sum_{j=0}^{\infty} h^{\frac{j}{2}}
 {\cal R}^{\infty}_{2-\frac{j}{2}}(s,D_s,x,D_x)$$
where
\medskip

(i) $ {\cal R}^{\infty}_{2-\frac{j}{2}}(s, D_s, x, D_x) =
  {\cal R}^{\infty,2}_{2-\frac{j}{2}} {\cal R}^2 +
   {\cal R}^{\infty,1}_{2-\frac{j}{2}} {\cal R} + {\cal
R}^{\infty,o}_{2-\frac{j}{2}},$ with
${\cal R}^{\infty,k}_{2-\frac{j}{2}} \in C^{\infty}(S^1_L, {\cal
E}_{\epsilon}^{j -2k});$
\medskip

(ii) $ {\cal R}^{\infty}_{2-j}(s, D_s, x, D_x)|_o = {\cal
R}^{\infty,o}_{2-j}(s, x, D_x)|_o =
  f_{j}(I_1^e,\dots,I_p^e,I^h_1,\dots,I^h_r, I_1^{ch Re},\dots I_c^{ch Re},
I_1^{ch Im}, \dots, I_c^{ch Im})|_o$ for certain polynomials $f_j$
 of degree j+2 on $\R^n,$ i.e.
 $f_{j}(I_1^e,\dots,I_p^e,I^h_1,\dots,I^h_r, I_1^{ch Re},\dots I_c^{ch Re},
I_1^{ch Im}, \dots, I_c^{ch Im} ) \in
{\cal P}^{j+2}_{\cal I}$
\medskip

(iii) $ {\cal R}^{\infty}_{2-\frac{2k+1}{2}}
(s, D_s, x, D_x)|_o = {\cal R}^{\infty, o}_{2-\frac{2k+1}{2}}
(s,x, D_x)|_o =
0;$

(iv) The terms $I^e_j, I^h_j, I^{ch}_j$  are weightless under metric scalings and all
of the ${\cal R}$'s have weight -2. }
\medskip

\noindent{\bf Proof}:
\medskip

 As in the elliptic case [Z.1,Lemma 2.22],
 the operator $W_h$ will be constructed as the asymptotic product
$$W_h:=\mu(P_{\gamma})^* \circ \Pi_{k=1}^{\infty}
W_{h \frac{k}{2}} \circ \mu(P_{\gamma}) \leqno(3.2)$$
of weightless unitary $h$-pseudodifferential operators on $\R^n$, with
$$W_{h \frac{k}{2}}:= exp(ih^{\frac{k}{2}} Q_{\frac{k}{2}})\leqno(3.3)$$ 
and with $h^{\frac{k}{2}} Q_{\frac{k}{2}}
 \in h^{\frac{k}{2}}\C^{\infty}(S^1_L)\otimes {\cal E}^{k+2}$
of total order 1.
The product will converge, for each s, to a unitary operator in 
$\Psi_h^o( \R^n)$
(we refer to [Sj] for a discussion of such asymptotic products).

  We first construct a weightless $Q_{\frac{1}{2}}(s, x, D_x)
 \in C^{\infty}(S^1_L)\otimes{\cal E}^3_{\epsilon}$ such that
$$e^{-ih^{\frac{1}{2}}Q_{\frac{1}{2}}} 
{\cal R}_h e^{ih^{\frac{1}{2}}Q_{\frac{1}{2}}}|_o
=[- h^{-2}L^{-2} + 
2 h^{-1}L^{-1}{\cal R} + {\cal R}^{\half}_o + \dots]|_o\leqno(3.4a)$$
where the dots $\dots$ indicate higher powers in $h$.  
The
 operator $Q_{\frac{1}{2}}$ then must satisfy
the commutation relation
$$\{ [L^{-1}{\cal R},Q_{\frac{1}{2}}]+
 {\cal R}_{\frac{1}{2}} \}|_o= 0.\leqno(3.4b)$$ 
 To solve for $Q_{\half}$,
we conjugate back to the ${\cal D}_{2 - \frac{m}{2}}$'s of the
twisted model
 by $\mu(P_{\gamma})$, which transforms ${\cal R}$
into $D_s$.  The commutation relation thus becomes
$$\{ [L^{-1} D_s,\mu(P_{\gamma}))^*Q_{\frac{1}{2}}\mu(P_{\gamma})]+ 
{\cal D}_{\frac{1}{2}} \}|_o= 0, \leqno(3.4c)$$
that is,
$$L^{-1}\partial_s\{\mu(P_{\gamma})^*Q_{\frac{1}{2}}\mu(P_{\gamma}))\}|_o =
 - i\{{\cal D}_{\frac{1}{2}}\}|_o\leqno(3.4d)$$
where $\partial_s A$ is the Weyl operator whose complete symbol is the $s$-derivative
of that of $A$.  Since (3.4d) is simpler than (3.4b), we henceforth conjugate
everything by $\mu(P_{\gamma}))$, and relabel the operators
$\mu(P_{\gamma})^*Q\mu(P_{\gamma})$ by
$\tilde{Q}.$  The resulting  ${\cal D}$'s  then transform under $\tau_L$ like
operators on the quantum mapping cylinder of $\mu(P_{\gamma})$.  Our problem is thus
to
 solve (3.4d) with an operator $\tilde{Q}_{\half}$ satisfying 
$$\tau_L \tilde{Q}_{\half} \tau_L^* = \mu(P_{\gamma}) \tilde{Q}_{\half}
\mu(P_{\gamma})^*.$$

To solve the equation (3.4d) we rewrite it in terms of  complete Weyl symbols.
We will use the notation $A(s,x,\xi)$ for the complete Weyl symbol of the
operator $A(s,x,D_x)$.  Then (3.4d) becomes
$$L^{-1}\partial_s \tilde{Q}_{\half}(s,x,\xi)= -i {\cal D}_{\half}|_o(s,x,\xi)
\leqno(3.5a)$$
with
$$\tilde{Q}_{\half}(s + L,x,\xi) = \tilde{Q}_{\half}(s,P_{\gamma}(x,\xi)).$$
We solve (3.5a) with the Weyl symbol
$$\tilde{Q}_{\half}(s,x,\xi) = \tilde{Q}_{\half}(0,x,\xi) + L \int_0^s
-i {\cal D}_{\half}|_o(u,x,\xi)du$$
where $\tilde{Q}_{\half}(0,x,\xi)$ is determined by the consistency condition
$$\tilde{Q}_{\half}(L,x,\xi) - \tilde{Q}_{\half}(0,x,\xi) = 
L \int_0^L -i {\cal D}_{\half}|_o(u,x,\xi)du \leqno(3.5b)$$
or in view of the periodicity condition in (3.5a),
$$\tilde{Q}_{\half}(0,P_{\gamma}(x,\xi)) - \tilde{Q}_{\half}(0,x,\xi) = 
L \int_0^L -i {\cal D}_{\half}|_o(u,x,\xi) du.\leqno(3.5c)$$

To solve, we use that ${\cal D}_{\half}|_o(u,x,\xi)$ is a polynomial of degree 3 in
$(x,\xi)$.  It will be most convenient to express this polynomial
in coordinates relative to the eigenvectors of the Poincare map.
In the elliptic planes, we use the complex
coordinates $z_j = x_j + i\xi_j$ and $\bar z_j = x_j - i \xi_j$ ($j=1, \dots, p$) in
which the action of $P_{\gamma}$ is diagonal. In the real hyperbolic planes we use the
real coordinates $(y_j, \eta_j) = (x_j, \xi_j), (j=p+1, \dots p+q)$ in which the real
hyperbolic part of
$P_{\gamma}$ is diagonal.  Finally, in the complex hyperbolic (loxodromic) 4-spaces
we use the coordinates $w_j = x_j + i x_{c + j}, \bar{w}_j = x_j -i x_{c+j}, \omega_j
=
\xi_j -i \xi_{c+j}, \bar{\omega_j} = \xi_j + i \xi_{c+j}, (j=p + q + 1, \dots p + q
+ c)$ in which the complex hyperbolic part of $P_{\gamma}$ is diagonal.

We will denote the Weyl symbols in these coordinates by their previous expressions.
 We also suppress the
subscripts by using vector notation $z, \bar z, y, \eta, w,\bar w, \omega, \bar 
\omega$.
 Thus, (3.5c) becomes
$$\tilde{Q}_{\half}(0,e^{i\alpha}z, e^{-i\alpha}\bar z, 
e^{\lambda}y,e^{-\lambda}\eta, e^{\mu + i\nu}w, e^{\mu - i\nu}\bar{w},
e^{-\mu + i\nu}\omega, e^{-\mu - i\nu}\bar{\omega} ) -
\tilde{Q}_{\half}(0,z,\bar z,y, \eta, w, \bar w, \omega, \bar \omega) =$$ 
$$=L\int_0^L -i {\cal D}_{\half}|_o(u,z,\bar z, y, \eta,w, \bar w, \omega, 
\bar \omega
) du.\leqno(3.6)$$ We now use that ${\cal D}_{\half}(u,z,\bar z, y, \eta,w, \bar w,
\omega, \bar \omega )$ is a polynomial of degree 3 to solve (3.5c).  If we put
$$\tilde{Q}_{\half}(s,z,\bar z,y,\eta,w, \bar w, \omega, \bar \omega) :=
\sum_{|a|+|\bar a|+|b| +|c| + |\bar c|\leq 3}
 q_{\half;a \bar a b c \bar c}(s) z^a\bar z^{\bar a} y^{b_1} \eta^{b_2}w^{c_1}
\omega^{c_2} \bar w^{\bar c_1} \bar{\omega}^{\bar c_2}\leqno(3.7a)$$
and 
$$ {\cal D}_{\half}|_o(s,z,\bar z,y,\eta,w, \bar w, \omega, \bar \omega ) du :=
\sum_{|a|+|\bar a|+|b|+|c|+ |\bar c|\leq 3} d_{\half;a \bar a b c \bar c}(s) z^a \bar
z^{\bar a} y^{b_1}
\eta^{b_2} w^{c_1}
\omega^{c_2} \bar w^{\bar c_1} \bar{\omega}^{\bar c_2}\leqno(3.7b)$$
then (3.6) becomes
$$\sum_{|a|+|\bar a|+|b|+|c|+ |\bar c|\leq 3} (1 - e^{i(a - \bar a) \alpha +
i(c_1 -\bar c_1) \nu +
(b_1 - b_2)\lambda + ( c_2 - \bar c_2)\mu }) q_{\half; a \bar a b c \bar c }(0)
z^a
\bar z^{\bar a} y^{b_1} \eta^{b_2} w^{c_1}
\omega^{c_2} \bar w^{\bar c_1} \bar{\omega}^{\bar c_2}=$$
$$ = -i L^2 \sum_{|a|+|\bar a|+|b|+|c|+ |\bar c|\leq 3} \bar d_{\half;a \bar a b c
\bar c} z^a \bar z^{\bar a} y^{b_1} \eta^{b_2} w^{c_1}
\omega^{c_2} \bar w^{\bar c_1} \bar{\omega}^{\bar c_2} \leqno(3.8).$$
Under the non-degeneracy assumption on $P_{\gamma}$, we can solve with
$$q_{\half; a \bar a b c \bar c}(0) = -i L^2 (1 - e^{i (a - \bar a)\alpha +
i(c_1 - \bar c_1) \nu +
(b_1 - b_2)\lambda + ( c_2 - \bar c_2)\mu})^{-1} d_{\half;a \bar a b c \bar
c}\leqno(3.9)$$ since $i(a - \bar a)\alpha +
i(c_1 - c_2) \nu )  + (b_1 - b_2)\lambda + (\bar c_1 - \bar c_2)\mu) = 2 \pi i
k$  only if $a=\bar a,b_1 = b_2, c= \bar c$ 
and there are no such $(a,\bar a, b_1, b_2,c,\bar c)$ in an odd-index equation.

Precisely as in the purely elliptic case of [Z.1], we see that
$\tilde{Q}_{\frac{1}{2}}$
is a pseudodifferential operator on
$\R^n$ with the same order,  same order of vanishing, and same parity 
as the restriction of ${\cal D}_
{\frac{1}{2}}$ to elements of weight zero.  We then extend it 
as a pseudodifferential operator of the form
$$\tilde{Q}_{\half} \in \Psi^o(\R^1) \otimes {\cal E}^3_{\epsilon}$$
 on all of  ${\cal H}_{\gamma}$ by decreeing that
it commute with $s$.
The conjugate by $\mu(P_{\gamma})$ then defines a unitary operator $W_{h \half} \in
\Psi^o_h(S^1 \times \R^n)$
 satisfying (3.4a).
The corresponding twisted unitary operator with exponent $\tilde{Q}_{\half}$,
i.e. the image of $W_{h \half}$ under
conjugation by
$\mu(P_{\gamma})$,  will be denoted
$\tilde{W}_{h \half}.$ 

The effect of this first conjugation is precisely as in the elliptic case:
 Since $h^{\half}\tilde{ Q}_{\half}$ is of total order 1,
     $h^{\half} ad(\tilde{ Q}_{\half})$ (with $ad(A)B:=[B,A]$)
 preserves the
total order in $\Psi_h^{(*,*,*)}$, and hence $\tilde{W}_{h \half}$ 
is an order-preserving automorphism of the model pseudodifferential algebra.
 It is moreover independent of $D_s$ and has an odd polynomial Weyl symbol, so that
$$h^{\half} ad(\tilde{Q}_{\half}): h^{\frac{k}{2}}\Psi^l(\R) \otimes {\cal E}^m_{\epsilon}
\rightarrow h^{\frac{k+1}{2}} [\Psi^{l-1}(\R) \otimes {\cal E}^{m+3}_{\epsilon} + 
\Psi^{l}(\R)
\otimes {\cal E}^{m+1}_{\epsilon}]. \leqno(3.10)$$
Finally, the $d_{\half;m,n}$'s have weight -2, the variables $z$ have weight 0 and
hence the $q_{\half;m,n}$'s have weight 0.

 Consider now the element
$${\cal D}^{\half}_h:=\tilde{W}_{h \frac{1}{2}}^*
 {\cal D}_h \tilde{W}_{h \frac{1}{2}} \in \Psi^2_h(\R^1\times
\R^n)$$
which can be expanded in the semi-classical series
$${\cal D}^{\half}_h\sim \sum_{n=o}^{\infty} h^{-2 + \frac{n}{2}} \sum_{j+m=n}
 \frac{i^j}{j!}
(ad\tilde{Q}_{\half})^j {\cal D}_{2 - \frac{m}{2}}\leqno(3.11)$$
$$:= h^{-2}L^{-2} + h^{-1}L^{-1}D_s +\sum_{n=3}^{\infty}
 h^{-2 + \frac{n}{2}}{\cal D}^{\half}_{2 - \frac{n}{2}}.$$ 
An obvious induction as in the elliptic case [loc.cit.] gives that
$$ad(\tilde{Q}_{\half})^j {\cal D}_{2 - \frac{m}{2}} \in C^{\infty}(\R,
 {\cal E}_{\epsilon}^{m+j-4})D_s^2
+ C^{\infty}(\R, {\cal E}_{\epsilon}^{m+j-2})D_s + C^{\infty}(\R, 
{\cal E}_{\epsilon}^{m+j}).$$
It follows that
${\cal D}^{\half}_{2 - \frac{n}{2}}$ has the same filtered structure  as 
${\cal D}_{2 - \frac{n}{2}}.$

We carry this procedure out one more step before referring to [Z.1] for the
inductive argument, since the
even steps behave differently from the odd ones.  
We thus seek an element $\tilde{Q}_1(s,x,D_x)\in \Psi^*(S^1 \times \R^n)$ and
 an element
$f_o(I_1^e,\dots,I_p^e, I_1^h, \dots, I_r^h, I^{ch, Re}_1,I^{ch,Im}_1,\dots,
I^{ch,Re}_c, I^{ch,Im}) \in
{\cal A}$ so that
$${\cal D}^1_h:=\tilde{W}_{h 1}^*{\cal D}^{\half}
 \tilde{W}_{h 1} = h^{-2}L^{-2} + h^{-1}L^{-1}D_s +
h^{-\half}{\cal D}^{\half}_{\half} + {\cal D}^1_o(s,D_s, x,D_x) + \dots
\leqno(3.12)$$ with 
$${\cal D}_o^1(s,D_s,x, D_x)|_o = f_o(I_1^e, \dots,
I_c^{ch,Im})\leqno(3.13a)$$ with $ \tilde{W}_{h 1}= e^{i h \tilde{Q}_1},$ and where
the dots signify terms of higher order in $h$. 
 Note that ${\cal D}^1_{\half} =
{\cal D}^{\half}_{\half}$, so that (3.12) implies that
$$\{h^{\half}{\cal D}^1_{\half} +{\cal D}_o^1\}|_o = f_o(I_1^e,\dots, I_c^{ch, Im}).
\leqno(3.13b)$$
 The condition on $\tilde{Q}_1$ is then
$$\{[D_s,\tilde{Q}_1] + {\cal D}_o^{\half}\}|_o = f_o(I_1^e,\dots, I_c^{ch, Im})
\leqno(3.14a)$$
 or equivalently
$$\partial_s\tilde{Q}_1|_o =\{-{\cal D}_o^{\half} + f_o(I_1^e, \dots,
I_c^{ch,Im})\}|_o \leqno(3.14b).$$

We solve (3.14b) by again
expressing everything in terms of complete Weyl symbols relative to
the eigenvector coordinates.  
  Thus we rewrite
(3.14b)) in the form
$$L^{-1} \partial_s \tilde{Q}_1(s,z, \bar z, y,\eta,w, \omega, \bar w,\bar \omega)
 =$$
$$=-i \{{\cal D}^{\half}_o|_o (s,z,\bar z, y,\eta,w, \omega, \bar w,\bar \omega) -
f_o(|z_1|^2,\dots, |z_p|^2, y_1\eta_1,\dots, y_r\eta_r,Re w\omega, Im w \omega)\}
\leqno (3.15a)$$ or equivalently
$$\tilde{Q}_1(s,z,\bar z, y,\eta,w, \omega, \bar w,\bar \omega) 
=\tilde{Q}_1(0,z,\bar z, y,\eta,w, \omega, \bar w,\bar \omega) $$
$$ -i L\int_0^s
[{\cal D}^{\half}_o|_o (u,z,\bar z, y,\eta,w, \omega, \bar w,\bar \omega) -
f_o(|z_1|^2,\dots, |z_p|^2,y_1\eta_1,\dots, y_r\eta_r,Re w\omega, Im w \omega)]
du\leqno (3.15b)$$   and solve simeltaneously for $\tilde{Q}_1$ and $f_o$. The
consistency condition determining a unique solution is that
$$\tilde{Q}_1(L,z,\bar z, y,\eta,w,\omega,\bar w,\bar \omega) =
 \tilde{Q}_1(0,z,\bar
z, y,\eta,w,\omega, \bar w, \bar \omega)$$
 $$-i  L\int_0^L [{\cal D}^{\half}_o|_o
(u,z,\bar z, y,\eta,w,\omega, \bar w, \bar \omega) - f_o(|z_1|^2,\dots,
|z_p|^2,y_1\eta_1,\dots, y_q\eta_q, Re w\omega, Im w \omega)] du.
\leqno(3.16a)$$ or in view of the twisted periodicity condition
$$\tilde{Q}_1(0,e^{i\alpha}z,e^{-i\alpha}\bar z,e^{\lambda}y,
e^{-\lambda}\eta,e^{\mu + i\nu} w,e^{\mu - i\nu}\omega, e^{- \mu + i\nu}\bar
w,e^{-\mu - i\nu}\bar \omega) -
\tilde{Q}_1(0,z,\bar z, y,\eta,w,\omega, \bar w, \bar \omega) =$$
$$ -i L 
\{\int_0^L {\cal D}^{\half}_o|_o (u,z,\bar z, y,\eta,w,\omega, \bar w, \bar
\omega)du  - L f_o(|z_j|^2,y_j\eta_j, 
Re w_m\omega_m, Im w_n \omega_n) \}.\leqno
(3.16b)$$  In the spirit of the previous step, we use that ${\cal
D}^{\half}_o|_o (u,z,\bar z, y,\eta,w,\omega, \bar w, \bar \omega)$
 is a polynomial of degree 4 to solve the equation.  We put
$$ \tilde{Q}_1(s,z,\bar z, y,\eta,w,\omega, \bar w, \bar \omega) =
\sum_{|a|+|\bar a| + |b| + |c| + |\bar c| \leq 4} q_{1;a\bar{a}bc\bar{c}}(s) z^a \bar
z^{\bar a} y^{b_1}\eta^{b_2}w^{c_1}\omega^{c_2}\bar{w}^{\bar c_1}
\bar{\omega}^{\bar c_2},$$
and in an abbreviated notation,
$$ f_o(|z|^2, y\cdot \eta, Re w \cdot \omega, Im w \cdot \omega) =
\sum_{|k|+|\ell| + |n| \leq 2}
 c_{o
k\ell n} |z|^{2k}(y\cdot \eta)^{\ell} (Re w\cdot \omega)^{n_1} (Im w\cdot
\omega)^{n_2},
\leqno(3.17)$$ and  
$${\cal D}^{\half}_o|_o (s,z,\bar z, y,\eta)du :=\sum_{|a|+|\bar a| + |b| + |c| +
|\bar c|\leq 4} d_{o; a\bar{a}bc\bar{c}}^{\half}(s) z^a \bar z^{\bar a}y^{b_1}
\eta^{b_2}w^{c_1}\omega^{c_2}\bar{w}^{\bar c_1}
\bar{\omega}^{\bar c_2},$$
and finally
$$\bar
d^{\half}_{o;a\bar{a}bc\bar{c}}:=
\frac{1}{L}\int_o^Ld^{\half}_{o;a\bar{a}bc\bar{c}}(s)ds $$ As above, we can solve for
the off-diagonal coefficients where either $a\not= b$ or $m\not= n$
$$q_{1; a\bar{a}bc\bar{c}}(0) = -i L^2 (1 - e^{i (a- \bar{a})\alpha+ i(c_1 -
\bar{c}_1) \nu + (b_1-b_2)\lambda) + (c_2 - \bar{c}_2)\mu})^{-1}
\bar d_{o; a\bar{a}bc\bar{c}}^{\half}
\leqno(3.18a)$$ and must set the diagonal coefficients with $a=\bar{a},c_1 =
\bar{c}_1, b_1 = b_2, c_2 = \bar{c}_2$ equal to zero. The expression in (3.18a)
is well-defined by the non-degeneracy assumption.
 The coefficients $c_{o k \ell}$
are then determined by
$$c_{o k \ell n} =   \bar d_{1; kk \ell \ell nn}^{\half}. \leqno(3.18b)$$

It is evident that $\tilde{Q}_1$ and $f_o(I_1^e,\dots,I_c^{ch,Im})$ are
even polynomial pseudodifferential  operators of degree 4 in the variables $(x,D_x)$,
 that $\tilde{Q}_1$
is weightless under metric rescalings and that the coefficients $c_{o k \ell n}$
are of weight -2.

The rest proceeds as in the elliptic case.\qed
\medskip

\section{ Normal form of the Laplacian: Proof of Theorem I}

We now use the semi-classical normal forms to put the Laplacian into quantum
Birkhoff normal form.  Essentially this amounts to taking direct sums (or integrals)
of the semi-classical normal form over various internal Planck constants.  
\medskip

\noindent{\bf Proof of Theorem I:}  As in the elliptic case, we make the transition
from the semi-classical normal form to the quantum Birkhoff normal form by using
genreralized eigenfunction expansions for the model algebra.  

  From the table in
\S 1.2.1 we
see that a function
$f \in L^2(S^1_L \times \R^p_x \times \R^q_y \times \R^{2c}_{r,\theta})$ can be
expanding in terms of joint ${\cal A}_{pqc}$-eigenfunctions as:
$$f(s,x,y,r, \theta) = \sum_{\pm} \sum_{(k,m, n) \in {\bf N}^{1 + p + c}}
\int_{\R^q}\int_{\R^{+ c}}
\hat{f}_{\pm} (k,m,n,a,t) e^{ir_{kmnat}} e^{i<n,\theta>}\gamma_m (x)y_{\pm}^{ia
-\half}r^{it-1} da dt.$$  Here as in \S 3, we have used the notation $x$ for
 linear coordinates on the elliptic factors, $y$ for those on the real hyperbolic
factors, and   polar coordinates $w_j = r_j e^{i
\theta_j}$ in each $w_j$-plane of the complex hyperbolic factors. We also employ
a multi-index notation.  

We now  assemble the semi-classical intertwining operators into the
Fourier-Hermite-Mellin -series-integral intertwining operator
$$W_{\gamma} : L^2(S^1_L \times \R^p\times \R^q \times \R^{2c}, ds dx dy dw),
\rightarrow L^2(S^1_L \times \R^p\times \R^q \times \R^{2c}, ds dx dy
dw)\leqno(4.1)$$
$$W_{\gamma} \sum_{\pm}\sum_{(k,m,n ) \in {\bf N}^{1+ p+c}}\int_{\R^q}\int_{\R^{+ c}}
\hat{f}(k,m,n,a,t) e^{ir_{kmnat}s}e^{i<n,\theta>}\gamma_m(x)y_{\pm}^{ia-\half}
r^{it-1}dadt =$$
$$\sum_{\pm}\sum_{(k,m,n) \in {\bf N}^{1+p+c}}\int_{\R^q}\int_{\R^{+ c}}
\hat{f}(k,m,n,a,t) e^{ir_{kmnrt}s}
W_{kmnrt}e^{i<n,\theta>}\gamma_m(x)y_{\pm}^{ir-\half}r^{it-1}dadt$$  with
$$W_{kmnrt}:= \mu(\tilde{{\cal W}}(s)^*) W_{r_{kqnat}^{-1}} \mu
(\tilde{{\cal W}}_s)^{*-1}.$$

 Also, the
dilation operators will be assembled into the  operator
$$T : L^2(S^1_L \times \R^p \times \R^q \times \R^{2c}, dsdxdydw) \rightarrow
L^2(S^1_L
\times \R^p \times \R^q \times \R^{2c}, dsdxdydw)),\leqno(4.2)$$ 
$$T \sum_{\pm}\sum_{(k,m,n) \in {\bf N}^{1+ p + c}}\int_{\R^q}\int_{\R^{+ c}}
\hat{f}(k,m,n,a,t) e^{ir_{kmnat}s}
e^{i<n,\theta>}\gamma_m(x)y_{\pm}^{ia-\half}r^{it-1}dadt=$$
$$ = \sum_{\pm}
\sum_{(k,m,n) \in {\bf N}^{1+p+c}}\int_{\R^q}\int_{\R^{+ c}} \hat{f}(k,m,n,a,t)
r_{kmnrt}^{i \half(|a|
+|t|)}e^{ir_{kmnat}s}e^{i<n,\theta>}\gamma_m(\sqrt{r_{kmnrt}}x)y_{\pm}^{ia-\half}
r^{it-1}dadt .$$  
Here, we used that the hyperbolic eigenfunctions are eigenfunctions of dilation
operators.

It follows, formally, from the semi-classical normal form and from the eigenfunction
expansion that
$$W_{\gamma}^{-1} T^{-1} \Delta T W_{\gamma} \sim
 {\cal L}^2 + f_o(I_{1}^e,...,I_{2c}^{ch, Im}) 
+ \frac{f_1(I_{1}^e,...,I_{2c}^{ch, Im})}{{\cal L}} + \dots.\leqno(4.3)$$
We now show
that the intertwining operator is actually a standard Fourier Integral operator
(in the Weyl operator, or isotropic, sense) and that (4.3) holds modulo the the
kind of error stated in Theorem I.

The proof is again similar to the elliptic case, so we  concentrate on the
novel aspects and refer the reader to [Z.1, Proposition 3.4] for the remaining
details. As before, we will not be as careful here as in [Z.1] to
express things in weightless terms relative to metric rescalings. 
\medskip

\noindent(4.4) \; {\bf Propostion} {\it $\; T W_{\gamma}T^{-1}$ is a
 (standard) Fourier integral operator, well-defined and invertible 
on the microlocal neighborhood (0.1) in $T^*(S^1_L \times \R^n)$.}
\medskip

\noindent{\bf Sketch of Proof}:  

We first consider the unitarily equivalent operator
$\tilde{T}\tilde{W}\tilde{T}^{-1}$ in the microlocal neighborhood (1.2.1) in the
twisted model, with
$$\tilde{W}: {\cal H}_{\alpha} \rightarrow {\cal H}_{\alpha} \leqno(4.5)$$
$$\tilde{W}(e^{ir_{kmnat}s}e^{i<n,\theta>}\gamma_m(x)y_{\pm}^{ia -\half} \rho^{it-1})
:= e^{ir_{kmnat}s}W_{r_{kmnat}^{-1}}e^{i<n,\theta>}\gamma_m(x)y_{\pm}^{ia-\half}
r^{it-1},$$ and with
$\tilde{T}$ the dilation operator analogous to (4.2) but relative to the basis
$e^{ir_{kmnat}s}e^{i<n,\theta>}\gamma_q(x)y_{\pm}^{ia-\half} r^{it-1}$.  
We then factor $\tilde{T}\tilde{W}\tilde{T}^{-1}$ as the product 
 $\tilde{T}\tilde{W}\tilde{T}^{-1}
= j^*\tilde{T} V \tilde{T}^{-1} $ where:
 $$ V: {\cal H}_{\alpha} \rightarrow L^2_{loc}(\R \times \R \times
\R^n)\leqno(4.6)$$
$$V=\Pi_{j=o}^{\infty} exp[i D_s^{-\frac{j}{2}} Q_{\frac{j}{2}}(s',y,D_y)]$$ that is,
$$Ve^{ir_{kmnat}s}e^{i<n,\theta>}\gamma_m(x)y_{\pm}^{ia-\half}
r^{it-1}):=e^{ir_{kmnat}s}W_{r_{kmnat}^{-1}}(s',x,y,w,
D_x,D_y,D_w)e^{i<n,\theta>}\gamma_m(x)y_{\pm}^{ia-\half}
r^{it-1},$$ and where
 $$j^*: C^{\infty}(\R \times \R \times \R^n) \rightarrow C^{\infty}(\R \times \R^n)
\leqno(4.7)$$
$$j^*f(s,x)=f(s,s,x)$$
is the pullback under the partial diagonal embedding. 

 The discussion of $V$ and the proof that $j^* V$ is a standard Fourier Integral
operators goes precisely as in [Z.1, Proposition 3.4].  The effect of the dilation
is to convert the isotropic calculus into the pure polyhomogenous calculus (see
also [G.1] for this aspect) and then the power of $D_s$ insures that the phases are
all homogeneous of degree 1 and vanishing to higher and higher order along $\gamma$
(by one step as the index j increases by one unit).  Hence the phase of the infinite
product has only finitely many terms of a given vanishing order and converges as
a formal power series in the transverse variable.  A convergent product can be
defined (by Borel summation) of the phase (cf. [Sj]).

The proposition then follows by expressing
$$TW_{\gamma}T^{-1} = T \mu({\cal W})^*\tilde{T}^{-1}
 \tilde{T} W \tilde{T}^{-1}\tilde{T}\mu({\cal W})T^{-1}$$
and noting that $\tilde{T}\mu({\cal W})T^{-1}$ is also
 a standard Fourier Integral operator. \qed

We now complete the proof of the quantum normal form Theorem I for $\sqrt{\Delta}$,
stated in an equivalent form in terms of $W_{\gamma}.$  
 As in the introduction,
the notation $A\equiv B$ means that the complete (Weyl) symbol of $A-B$ vanishes to
infinite order at $\gamma$ and  $O_j\Psi^m$ denotes the pseudodifferential
operators of order m whose Weyl symbols vanish to order j at $(y,\eta)=(0,0)$.  Here,
pseudodifferential operator can refer to either the standard polyhomogeneous kind, or
to the mixed polyhomogeneous-isotropic kind as in $\Psi^k(S_L^1)\otimes {\cal W}^l$,
in which case the total order is defined to be $m=k+l$. To simplify notation, we will
denote the space of mixed operators of order m by $\Psi_{mx}^m(S^1_L\times \R^n)$.
\medskip
 
\noindent(4.8) \;{\bf Lemma} {\it $\;\;$  Let $TW_{\gamma}T^{-1}$ be the Fourier
Integral operator of Proposition (4.4), 
defined over a conic neighborhood of $R^+\gamma$ in $T^*(S^1_L \times \R^n)$.  Then:
   $$W_{\gamma}^{-1} T^{-1} \sqrt{\Delta} T W_{\gamma} \equiv
 P_{1}({\cal L},I_{1}^e,...,I_{2c}^{ch, Im}) 
+ P_{o}({\cal L},I_{1}^e,...,I_{2c}^{ch, Im}) + \dots \;\;\;\mbox{mod}\;\;\;\;
\oplus_{k=o}^{m+1}
O_{2(m+1-k)}\Psi^{1-k}_{mx}(S^1_L\times \R^n),$$
where
$$P_1({\cal L}, I_{1}^e,...,I_{2c}^{ch, Im})\equiv {\cal L} +
\frac{ p_{1}^{[2]}(I_{1}^e,...,I_{2c}^{ch, Im})}{L {\cal L}} 
+ \frac{p_{2}^{[3]}(I_{1}^e,...,I_{2c}^{ch, Im})}{(L{\cal L})^2} + \dots
\leqno(4.9)$$
$$P_{-m}({\cal L}, I_{1}^e,...,I_{2c}^{ch, Im})\equiv \sum_{k=m}^{\infty} 
 \frac{p_{k}^{[k-m]}(I_{1}^e,...,I_{2c}^{ch, Im})}{(L{\cal L})^j}$$
with $p_{k}^{[k-m]}$, for m=-1,0,1,..., homogenous of degree l-m in the variables
$(I_{1}^e,...,I_{2c}^{ch, Im}) $ and of weight -1.  }
\medskip

\noindent{\bf Proof}:  

 As a semi-classical expansion in the ``parameter" $h = \frac{1}{L {\cal L}}$,
 (4.3) may be rewritten in the form :
$$W_{\gamma}^{-1} T^{-1} \sqrt{\Delta} T W_{\gamma} \sim
 {\cal L} +\frac{ p_{1}(I_{1}^e,...,I_{2c}^{ch, Im})}{L{\cal L}} 
+ \frac{p_{2}(I_{1}^e,...,I_{2c}^{ch, Im})}{(L{\cal L})^2} + \dots.\leqno(4.10)$$ 
From the fact that the numerators
$f_j(I_{1}^e,...,I_{2c}^{ch, Im})$ in (4.3) are polynomials of degree j+2 
and of weight -2, the
numerators $ p_{k}(I_{1}^e,...,I_{2c}^{ch, Im})$ are easily seen to be
 polynomials of degree $k+1$ and of weight -1. Hence they may be expanded 
in homogeneous terms 
$$p_k = p_k^{[k+1]} + p_k^{[k]} + \dots p_k^{[o]}, \leqno(4.11)$$
 with  $p_k^{[j]}$ the term of degree j and still of weight -1.
  The right side of (4.12) can
then be expressed as a sum of homogeneous operators:
$$P_1({\cal L}, I_{1}^e,...,I_{2c}^{ch, Im}) + P_o({\cal L},
I_{1}^e,...,I_{2c}^{ch, Im})
+ \dots\leqno(4.12)$$
with
$$P_1({\cal L}, I_{1}^e,...,I_{2c}^{ch, Im})\equiv {\cal L} +
\frac{ p_{1}^{[2]}(I_{1}^e,...,I_{2c}^{ch, Im})}{L{\cal L}} 
+ \frac{p_{2}^{[3]}(I_{1}^e,...,I_{2c}^{ch, Im})}{(L{\cal L})^2} + \dots
\leqno(4.13)$$
$$P_{-m}({\cal L}, I_{1}^e,...,I_{2c}^{ch, Im})\equiv \sum_{k=m}^{\infty} 
 \frac{p_{k}^{[k-m]}(I_{1}^e,...,I_{2c}^{ch, Im})}{(L{\cal L})^k}.$$

We claim that:
$$W_{\gamma}^{-1} T^{-1} \sqrt{\Delta} T W_{\gamma} -
[ {\cal L} +\frac{ p_{1}(I_{1}^e,...,I_{2c}^{ch, Im})}{L{\cal L}} 
+ \frac{p_{2}(I_{1}^e,...,I_{2c}^{ch, Im})}{(L{\cal L})^2} +\dots
+\frac{p_{m}(I_{1}^e,...,I_{2c}^{ch, Im})}{(L{\cal L})^m}]\leqno(4.14)$$
$$ \in \oplus_{k=o}^{m+1}
O_{2(m+1-k)}\Psi^{1-k}_{mx}(S^1_L\times \R^n).$$

Indeed, from the analysis of the remainder terms in the semi-classical normal
form (see Lemma (3.1 (i)) and [Z.1, Lemma 2.22]), we have
$$P_1({\cal L}, I_{1}^e,...,I_{2c}^{ch, Im}) - [{\cal L} +
\frac{ p_{1}^{[2]}(I_{1}^e,...,I_{2c}^{ch, Im})}{L{\cal L}} 
+ \dots + \frac{p_{N}^{[N+1]}(I_{1}^e,...,I_{2c}^{ch, Im})}{(L{\cal
L})^k}]\leqno(4.15)$$
$$ \in O_{2(N+2)}\Psi^1_{mx}(S^1_L\times \R^n)$$
and also

$$P_{-m}({\cal L}, I_{1}^e,...,I_{2c}^{ch, Im}) - \sum_{k=m}^{N} 
 \frac{p_{k}^{[k-m]}(I_{1}^e,...,I_{2c}^{ch, Im})}{(L{\cal L})^k}\in O_{2(N+1-m)}
\Psi^{-m}_{mx}(S^1_L\times \R^n).\leqno(4.16)$$

Hence the expansion (4.10) is also asymptotic in the sense of $\equiv.$  For the
statement of Theorem I in the introduction, it is only necessary to conjugate under
$\mu(\tilde{W}).$ The rest proceeds as in the elliptic case.\qed
\medskip

\section{Wave invariants and residue trace: Proof of Theorem B}

The purpose of this section is to show that the wave invariants have precisely the
same relation to the coefficients of the quantum normal form in the non-degenerate
case that they have in the elliptic one.  The characterization of the wave invariants
in Theorem I will then follow from Theorem A of [Z.1].

We will need to use some further notation and results from [Z.1]: First, the
kth wave invariant of a positive elliptic operator $P$ at a non-degenerate
closed bicharacteristic $\gamma$ will be denoted $\tau_{\gamma k}(P).$  According
to [Z.1, Proposition 4.2] we then have:
$$\tau_{\gamma k}(P) = \tau_{\gamma k}(P_1^{\leq 2k + 4} + P_o^{\leq 2k + 2}  + 
\dots + P_{-k-1}^{o}) \leqno(5.1) $$
where $P_j^{\leq k}$ denotes the first k terms in the Taylor expansion of the
jth homogeneous part of the complete symbol of $P$ at $\gamma$. Thus,   
$\tau_{\gamma k}(P)$ involves the (2k+4)th jet of the principal symbol, the
(2k+2)-jet
of the subprincipal term, ..., up to the zero-jet of term of
 homogeneity order  (-k-1).

As in [Z.1, \S 4],  we will also  rewrite the normal form in terms of
 of $D_s$ and $H_{\alpha,\lambda, (\mu,\nu)}$ using that
$$\frac{p_{\nu}(I_1^e,\dots,I_{2c}^{ch, Im})}{(L{\cal R})^{\nu}} =
\frac{p_{\nu}(I_1^e,\dots,I_{2c}^{ch, Im})} {(LD_s)^{\nu}} (I -\nu
\frac{H_{\alpha,\lambda, (\mu,\nu)}}{L D_s} +\half
\nu (\nu-1)  (\frac{H_{\alpha,\lambda, (\mu,\nu)}}{L D_s})^2 + \dots).$$
By (5.1), we can drop the  $D_s^{-(k+1) + \nu}H_{\alpha,\lambda, (\mu,\nu)}^{k+1
-\nu}$ and higher terms, so ${\cal D}_{k+1}$ can be written in the
form 
$${\cal D}_{k+1} \equiv L D_s + H_{\alpha,\lambda, (\mu,\nu)}+
\frac{\tilde{p}_1(I_1^e,\dots,I_{2c}^{ch, Im})}{L D_s} +
 \frac{\tilde{p}_2(I_1^e,\dots,I_{2c}^{ch, Im})}{(LD_s)^2}
+\dots+\frac{\tilde{p}_{k+1}(I_1^e,\dots,I_{2c}^{ch, Im})}{(LD_s)^{k+1}}
\leqno(5.2)$$ modulo terms which make no contribution to $\tau_{k\gamma}$.

We then use the fact ([Z.1, (4.3)],[Z.2]) that 
$$\tau_{\gamma k}(P) = res D_t^k \psi_{\epsilon}(D_s,y,D_y) e^{ it P}|_{t=L}
\leqno(5.3)$$  where $res$ is the non-commutative residue and where
$\psi_{\epsilon}(D_s,y,D_y)$ denotes a microlocal cut-off to the cone (1.2.1).
Note that in contrast to the elliptic case, the microlocal cut-off cannot be
constructed in ${\cal A}_{p,q,c}$ since the neighboorhoods given by $I < \epsilon
\sigma$ in terms of mixed hyperbolic-elliptic actions are  of infinite
transverse symplectic volume. This does not pose a genuine problem, but accounts
for a number of modifications to the elliptic case in [Z.1].
  For the gauging elliptic operator  we  use $L D_s$.

  To simplify the
notation we will put 
$${\cal P}_{k+1}(D_s, I_1^e,\dots,I_{2c}^{ch,
Im}):=\frac{\tilde{p}_1(I_1^e,\dots,I_{2c}^{ch, Im})}{L D_s} +
 \frac{\tilde{p}_2(I_1^e,\dots,I_{2c}^{ch, Im})}{(LD_s)^2}
+\dots+\frac{\tilde{p}_{k+1}(I_1^e,\dots,I_{2c}^{ch, Im})}{(LD_s)^{k+1}}.
\leqno(5.4)$$ 
so that:
$$\tau_{\gamma k}(\sqrt{\Delta}) = Res_{z=0}
 Tr D_t^k \psi_{\epsilon}(D_s,y,D_y)
e^{it[\frac{1}{L} (2 \pi L D_s + H_{\alpha,\lambda, (\mu,\nu)} )+ {\cal P}_{k+1}]}
(LD_s)^{-z} |_{t=L}.
\leqno(5.5)$$

As in the elliptic case, a key role will be played by the (formal) trace
$$T(\alpha,\lambda,(\mu,\nu)):= Tr e^{iH_{\alpha,\lambda, (\mu,\nu)}}.$$
Its precise definition is the following:
  Since $e^{iH_{\alpha,\lambda, (\mu,\nu)}}$ is
an element of the metaplectic representation $\mu$ of $Mp(n,\R)$, 
$T(\alpha,\lambda,(\mu,\nu))$ may be 
identified with the character $Ch$ of $\mu$ evaluated at the associated
element
$P_{\gamma} = exp( \Xi_{H_{\alpha,\lambda, (\mu,\nu)}}) \in Mp(n, \R).$  Here,   
$H_{\alpha,\lambda, (\mu,\nu)}$ denotes the quadratic function on $\R^{2n}$
which gives the complete Weyl symbol of the corresponding action operator, and
as above $exp \circ \Xi$ denotes the flow at time 1 of its Hamilton vector
field, or, more correctly, the lift to $Mp(n, \R)$  which
corresponds to $e^{iH_{\alpha,\lambda, (\mu,\nu)}}$ under $\mu$.

Since $Mp(n, \R)$ is a semi-simple Lie group,
the character $Ch$ is a real analytic function on the open dense
subset $Mp(n,\R)_{reg}$  of regular elements of $Mp(n, \R),$ where it is given
by the Harish-Chandra formula [Kn]. We will need below the explicit formula for
$Ch(x)$ in terms of the eigenvalues of $x$.  For elements of $Mp(n,\R)$ not having  1
as an eigenvalue, we recall that $Ch(x)$ is 
 given by
$$Ch (x) = \frac{i^{\sigma}}{\sqrt{|det (I - x)|}}$$
where $\sigma$ is a certain Maslov index. 
  For non-degenerate $x$ with p pairs of
eigenvalues $e^{\pm i \alpha_j}$ of modulus one, q pairs of positive real 
eigenvalues $e^{\pm \lambda_j}$ and c quadruplets of eigenvalues $e^{\pm  (\mu_j
\pm i \nu_j)}$, $Ch(x)$ is therefore given (up to a Maslov factor) by
$$T(\alpha,\lambda,(\mu,\nu)) = \Pi_{j=1}^p \frac{e^{\half i\alpha_j}}{1 -
e^{i\alpha_j}} \cdot \Pi_{j=1}^q \frac{e^{\half \lambda_j}}{1 -
e^{\lambda_j}} \cdot \Pi_{j=1}^{c} \frac{e^{\half(\mu_j + i\nu_j)}}{1 -
e^{\mu_j + i\nu_j}}\frac{e^{\half(\mu_j - i\nu_j)}}{1 -
e^{\mu_j - i\nu_j}}.\leqno(5.6)$$
Here we have selected one eigenvalue $\rho$ from each symplectic pair $\rho,
\rho^{-1}$ (see \S 1.1-2).  The ambiguity is fixed by the Maslov factor $i^{\sigma}$,
which can (and will) be ignored below for the sake of brevity.

 We can now give:
\medskip

\noindent{\bf Proof of Theorem B:}  Since $e^{2\pi i L D_s}
\equiv I$ on 
$L^2(S^1_L)$ we have
$$a_{k \gamma} = \tau_{\gamma k}(\sqrt{\Delta}) =$$
$$Res_{z=0} Tr \psi_{\epsilon}(D_s,y,D_y) [\frac{1}{L}2 (\pi L  D_s +
H_{\alpha,\lambda, (\mu,\nu)}) + {\cal P}_{k+1}]^k e^{i H_{\alpha,\lambda, (\mu,\nu)}}
 e^{iL{\cal P}_{k+1}} (LD_s)^{-z}. \leqno(5.7)$$ 
In view of the microlocal cutoff, 
 the operator under the trace is of trace class for $Re z$ sufficiently large.
Indeed, in estimating the trace we may eliminate the unitary factors and we are
then left with a pseudodifferential operator whose complete symbol is a polynomial
in $(\sigma, y,\eta)$ times a factor of $\sigma^{-Re z} \psi_{\epsilon}(\sigma,
y,\eta).$  The integral in the transverse $(y,\eta)$ variables is bounded by the
volume of the ball $(y^2 + \eta^2) < \sigma$ and hence is of order
$\sigma^{n}$.  Since a pseudodifferential operator is Hilbert-Schmidt if its Weyl
symbol is in $L^2$, the operator under the trace is Hilbert-Schmidt for $Re z >n$
and in particular is of trace class.  Moreover, since it is the non-commutative
residue of a Fourier Integral operator, one knows apriori that it admits a
meromorphic continuation to
$\C$ with at most simple poles [Z.2]. Hence the residue is well defined.

 As in [Z.1], we view the trace as a function
of the parameters $(\alpha, \lambda, \mu,\nu)$ and use the explicit form of the
exponential in $e^{i H_{\alpha,\lambda, (\mu,\nu)}}$ to rewrite (5.7) in the form
$$Res_{z=0}  \sum_{n=1} ^{\infty}  n^{-z} Tr \psi_{\epsilon}(D_s,y,D_y) \{
[\frac{1}{L} (2\pi n + \sum_{j=1}^p\alpha_j D_{\alpha_j} +
\sum_{j=1}^q\lambda_j \partial_{\lambda_j} +$$
$$+ \sum_{j=1}^{2c}(\mu_j
\partial_{\mu_j} + \nu_j D_{\nu_j} ) + {\cal P}_{k+1} (n, D_{\alpha_1},\dots,
D_{\nu_{2c}},L)]^k e^{iL{\cal P}_{k+1}(n, D_{\alpha_1},\dots, D_{\nu_{2c}},L)}
 e^{i H_{\alpha, \lambda, (\mu,\nu)}}\} .\leqno(5.8)$$
Here, we have used that $D_s$ commutes with $(y,D_y)$ to replace it by its eigenvalue
in the $s$-trace, and we have repeatedly used identities of the form
$F(D_x) e^{ix P} = F(P) e^{ix P}$ ($x\in \R$). 

Since ${\cal P}_{k+1}(n, D_{\alpha_1},\dots, D_{\nu_{2c}},L)$ is a symbol of order
$-1$ in $n$ with coefficients given by polynomials in the operators $D_{\alpha_j}$
(etc.),  we can expand the kth power in (5.8)   as an
operator-valued polyhomogeneous function of $n$.  At least formally, we can also
expand the exponential $e^{iL{\cal P}_{k+1}(n, D_{\alpha_1},\dots, D_{\nu_{2c}},L)}$
in a power series and then expand each term in the power series as a polynomial in
$n^{-1}$. 
 Collecting
powers of n, the right side of (5.8) may be written in the form
$$Res_{z=0} \sum_{n=1}^{\infty} \sum_{j= o}^{\infty} n^{-z+ k -j} {\cal F}_{k,k-j}
(D_{\alpha_1},\dots, D_{\mu_c + i \nu_c}, D_{\mu_c - i
\nu_c})
Tr \psi_{\epsilon}(n ,y,D_y) e^{i H_{\alpha, \lambda, (\mu,\nu)}}, 
\leqno(5.9)$$ with ${\cal F}_{k,k-j}(D_{\alpha},\dots, D_{\mu_c + i \nu_c}, D_{\mu_c - i
\nu_c}L)$ the coefficient of
 $n^{k-j}$
in (5.8).  The expansion of the exponential is justified
as in the elliptic case: as in [Z.1, (4.29)] we may write
$$e^{iL{\cal P}_{k+1}}:=  e_N(i L {\cal P}_{k+1}) +
(i L {\cal P}_{k+1})^{N+1} b_N(i L {\cal P}_{k+1})$$
with $e_N(ix) = 1 + ix + \dots + \frac{(ix)^N}{N!}$, with ${\cal P}_{k+1}$ short for 
${\cal P}_{k+1} (n, D_{\alpha_1}, \dots, D_{\nu_{2c}})$ and with $b_N(ix)$ a bounded
function.
The $e_N$ term  contributes a finite number of terms of
 the desired form (5.9).
For the remainder,
we expand  $(i L {\cal P}_{k+1})^{N+1}$ as a polynomial in $n^{-1}$ 
with coefficients given by operators $Q_{N p}(D_{\alpha_1},\dots,D_{\nu_{2c}})$ 
and observe that
each term has a factor of $n^{-N-1}$.  For each such term, we
 remove the coefficient operator
$Q_{N p}$
from the sum $\sum_n$, as above, leaving only the factor of $b_N$.  
  Since $b_N(ix)$ is a bounded function, it follows that $b_N(i L {\cal P}_{k+1})$ is
a bounded operator on $L^2$; and since each term of the resulting sum has at least the
factor 
$n^{-z -N -1 + k}$ (possibly multiplied by a further negative power of $n$), we see
that the remainder is a sum of terms of the form
$$Res_{z=0} Q_{N p}( D_{\alpha_1},\dots,D_{\nu_{2c}}) \sum_{kn} n^{-z +k -N -1 -l}
b_N(i{\cal P}_{k+1}(n, D_{\alpha_1},\dots))
Tr \psi_{\epsilon}(n ,y,D_y) e^{i H_{\alpha, \lambda, (\mu,\nu)}}.\leqno(5.10)$$

We then observe that the sum is bounded by $\sum_{m=1}^{\infty} m^{-Re z - N - 1 + k +n}$,
hence converges absolutely and uniformly for $Re z > -N + k + n$.  It follows that
for $N > (n+k)$ the sum in (5.10) defines a holomorphic function of $z$ in a
half-plane containing $z=0$ and since the operations of taking the residue in $z$
and  derivatives in $\alpha$ commute, each term (5.10) is zero.  This justifies
(5.10) and shows that 
 it is actually a finite sum in j, say $j < M$ (in fact M=(k+1)(n+k+1)).

The residue in (5.10) is therefore well-defined  and independent of $\epsilon$. Since
$Tr \psi_{\epsilon}(n ,y,D_y) e^{i H_{\alpha, \lambda, (\mu,\nu)}} \rightarrow
T(\alpha,\lambda, \mu,\nu)$ in the sense of distributions as
$\epsilon \rightarrow
\infty$  we must have
$$a_{\gamma k} =Res_{z=0} \sum_{m=1}^{\infty} \sum_{j= o}^{M} m^{-z + k -j} 
{\cal F}_{k,k-j}(D_{\alpha},\dots,D_{\mu_c + i \nu_c}, D_{\mu_c - i
\nu_c})\;\;
T(\alpha,\lambda, \mu,\nu) \leqno(5.11)$$  
$$=Res_{z=0} \sum_{j=0}^{M} \zeta(z + j - k){\cal F}_{k,k-j}(D_{\alpha},\dots,D_{\mu_c + i \nu_c}, D_{\mu_c - i
\nu_c})\;\;
T(\alpha,\lambda, \mu,\nu).$$ Here, $\zeta$ is the Riemann zeta-function, which has
only a simple pole at $s=1$ with reside equal to one.  It follows that the only term
contributing to (5.11) is that with
$j = k + 1$ and hence we have
$$a_{\gamma k} =  {\cal F}_{k, -1}(D_{\alpha},\dots, D_{\mu_c + i \nu_c}, D_{\mu_c - i
\nu_c}) T(\alpha,\lambda, \mu,\nu).
\leqno(5.12)$$
It follows that the wave invariants consist of the geometric data contained in
the coefficients of the ${\cal F}_{k, -1}$'s, and hence in the normal form
coefficients.  But the algorithm for constructing the normal form is essentially
the same as in the elliptic case and so the geometric data entering into the
normal form coefficients is of precisely the same kind.\qed

\section{Inverse Problems: Proofs of Theorem II and Corollary II.1}

Our first goal in this section is to prove that the wave invariants of 
$\gamma, \gamma^2,\dots$ determine the quantum normal form coefficients at
$\gamma.$
\medskip

\noindent{\bf Proof of Theorem II}:

We recall that the quantum normal coefficients are the coefficients of the
action monomials in the 
action polynomials $\tilde{p}_j(I_1^e, \dots, I_{2c}^{ch,Im})$ of (5.2).  These
coefficients determine, and are determined by, the
coefficients of the monomials in the polynomials $p_{\nu} (I_1^e, \dots,
I_{2c}^{ch,Im})$ in Theorem B.  They also corresponding bi-uniquely to the
coefficients of the classical action monomials in the complete symbols of either
set of action polynomials.

We also observe that the quantum normal form coefficients determine, and are
determined by, the coefficients of the constant coefficient partial differential
operator (PDO) 
$$ {\cal F}_{k,-1}(D_{\alpha}, D_{\lambda}, D_{\mu + i \nu},
D_{\mu - i
\nu}) := \sum_{(a, b,c_1,c_2) \in \Nb^{n}: |a| + |b| + |c_1| + |c_2| \leq k+1} 
C_{k; abc_1 c_2} D_{\alpha}^{a} D_{\lambda}^b D_{\mu + i \nu}^{c_1}  D_{\mu - i
\nu}^{c_2}
\leqno(6.1 k)$$ where $a \in \Nb^p, b\in \Nb^q, c_1,c_2 \in \Nb^c.$ 
 This can be proved easily by induction on k: In the case k=1, 
${\cal F}_{k,-1}$ is obtained from $\tilde{p}_{1}$ by substituting the variables
$D_{\alpha_1},$ etc. in for the variables $I^e_1,$ etc. Assuming inductively that
we have determined the coefficients of $\tilde{p}_1,\dots, \tilde{p}_k$ from those
of ${\cal F}_{1,-1},\dots, {\cal F}_{1,k}$, we note that $\tilde{p}_{k+1}$
contributes to the residue (5.7) for the first time at the $k+1$st stage. Since
it comes with the denominator $D_s^{k+1}$, it  only contributes to the residue
when composed with $D_s^k$.  From the form of (5.7) it is clear that only one
term involving $\tilde{p}_{k+1}$ contributes non-trivially, and that is the one
which appears in the linear term in the expansion of the exponential. Hence,
its contribution to ${\cal F}_{k+1,-1}$ is again just the substitution of
the variables $D_{\alpha_1},$ etc. in for the variables $I^e_1,$ etc.

  We next observe that the quantum normal form of $\Delta$ at any iterate
$\gamma^N$  of $\gamma$ is the same as for the primitive $\gamma$ itself.  Hence the
PDO ${\cal F}_{k,-1}$ is independent of the number $N$ of 
iterations. 
On the other hand, under the iteration $\gamma \rightarrow \gamma^N$, the Poincare map
transforms by $P_{\gamma^N}
\rightarrow P_{\gamma}^N$.  Therefore the expression
in (5.12) for the kth wave invariant of $\gamma^N$ is given by:
$${\cal F}_{k,-1}(D_{\alpha'}, D_{\lambda'}, D_{\mu' + i \nu'}, D_{\mu' - i
\nu'}) \cdot$$
$$\cdot \Pi_{j=1}^p \frac{e^{\half i  \alpha_j'}}{(1 - e^{i  \alpha_j'})} \cdot
\Pi_{j=1}^q
\frac{e^{\half  \lambda_j'}}{(1 - e^{  \lambda_j'})} \cdot \Pi_{j=1}^{c}
\frac{e^{\half  (\mu_j' + i\nu_j')}}{(1 - e^{ (\mu_j' + i\nu_j')})}\frac{e^{\half
(\mu_j' - i\nu_j')}}{(1 - e^{(\mu_j' - i\nu_j')})}|_{(\alpha', \lambda',
\mu', \nu') = N ( \alpha, \lambda, \mu, \nu)}.\leqno(6.2 k, N)$$ It therefore suffices
to prove  that for all k the coefficients of the  PDO ${\cal F}_{k,-1}$ can be
determined from its values (6.2 k, N) on $T(N\alpha, N\lambda, N \mu, N\nu)$ for $N =
\pm 1, \pm 2,\dots.$

We begin the proof by noting that (6.2 k, N) can be rewritten as:  
$$\Pi_{j=1}^p e^{\half i N \alpha_j}\cdot
\Pi_{j=1}^q e^{\half N \lambda_j}\cdot \Pi_{j=1}^{c}e^{\half N (\mu_j +
i\nu_j)} e^{\half N(\mu_j - i\nu_j)}\cdot \leqno(6.3 k,
N)$$
$$\cdot {\cal F}_{k,-1}(D_{\alpha'} +
\half , D_{\lambda'} + \half, D_{\mu' + i \nu'} +\half , D_{\mu' - i
\nu'} + \half )[\Pi_{j=1}^p \frac{1}{(1 - e^{i \alpha_j'})}$$
$$ \cdot \Pi_{j=1}^q
\frac{1}{(1 - e^{  \lambda_j'})} \cdot \Pi_{j=1}^{c}
\frac{1}{(1 - e^{ (\mu_j' + i\nu_j')})}\frac{1}{(1 - e^{(\mu_j' - i\nu_j')})}]
|_{(\alpha', \lambda',
\mu', \nu') = N ( \alpha, \lambda, \mu, \nu)}.$$
Making the substitutions $D_{\alpha} \rightarrow D_{\alpha} + \half$ (etc.)
in (6.1 k) we obtain
a new PDO whose coefficients $C'_{k; ab c_1 c_2}$  correspond in a bi-unique way
with the original $C_{k; ab c_1 c_2}$'s. Hence it will suffice to show that
we can determine the $C_{k; ab c_1 c_2}$'s from the values (6.3 k,N).

To do so, we will regard (6.3 k, N) as the values at integral points $z = N$ of
a function of $z$.
From the fact that
$$D_{\alpha} (1 - e^{i \alpha})^{-1} =  [ (1 - e^{i \alpha})^{-2} - (1 - e^{i
\alpha})^{-1}]$$
we see  that this function is a  polynomial in 
$(1 - e^{i z\alpha_j})^{-1}, (1 - e^{ z\lambda_j})^{-1}, (1 - e^{z(\mu + i
\nu}))^{-1}, (1 - e^{z(\mu - i \nu}))^{-1}.$ We clear the denominators to 
 obtain the entire function
$$\sum_{(a, b,c_1,c_2) \in \Nb^{n}: |a| + |b| + |c_1| + |c_2| \leq k+1} 
C'_{k; abc_1 c_2}[\Pi_{j=1}^p (1 - e^{i z \alpha_j}) \cdot \Pi_{j=1}^q\cdot
\leqno(6.4 k,z )$$
$$\cdot (1 - e^{ z \lambda_j}) \cdot \Pi_{j=1}^{c}
(1 - e^{z (\mu_j + i\nu_j)})(1 - e^{z (\mu_j - i\nu_j)})]^{(k+1)}
\cdot \Pi_{j=1}^p e^{\half i z \alpha_j}\cdot
\Pi_{j=1}^q e^{\half z \lambda_j}\cdot \Pi_{j=1}^{c}e^{\half z (\mu_j +
i\nu_j)} e^{\half z(\mu_j - i\nu_j)} \cdot$$ 
$$\cdot D_{\alpha'}^{a} D_{\lambda'}^b D_{\mu' + i \nu'}^{c_1}  D_{\mu' - i
\nu'}^{c_2}\cdot
\cdot [\Pi_{j=1}^p \frac{1}{(1 - e^{i  \alpha_j'})} \cdot
\cdot \Pi_{j=1}^q
\frac{1}{(1 - e^{  \lambda_j'})} \cdot \Pi_{j=1}^{c}
\frac{1}{(1 - e^{(\mu_j' + i\nu_j')})}\frac{1}{(1 - e^{(\mu_j' - i\nu_j')})}]
|_{(\alpha', \lambda',
\mu', \nu') = z ( \alpha, \lambda, \mu, \nu)}$$
which is an exponential polynomial  of the form
$$\sum_{\beta\in \Nb^{n}} c_{k; \beta} e^{ z
\langle
\beta + (\half,\dots,\half), (\alpha,
\lambda\mu, \nu) \rangle}. \leqno(6.5 k) $$ 
\medskip

\noindent{\bf (6.6) Lemma 1}~~~{\it The coefficients $C_{k; ab c_1 c_2}$  can be
determined from the coefficients $c_{k; \beta}$ in (6.5 k).}
\medskip

\noindent{\bf Proof}:  We first show that the coefficients $C_{k; a b c_1 c_2}$
with $|a| + |b| + |c_1| + |c_2| = k+1$ can be determined from the
$c_{k;\beta}$'s. Indeed, before multiplying by 
$$[\Pi_{j=1}^p (1 - e^{i z \alpha_j}) \cdot \Pi_{j=1}^q
 (1 - e^{ z \lambda_j}) \cdot \Pi_{j=1}^{c}
(1 - e^{z (\mu_j + i\nu_j)})(1 - e^{z(\mu_j - i\nu_j)})]^{(k+1)}\cdot$$
$$\cdot \Pi_{j=1}^p e^{\half i z \alpha_j}\cdot
\Pi_{j=1}^q e^{\half z \lambda_j}\cdot \Pi_{j=1}^{c}e^{\half z (\mu_j +
i\nu_j)} e^{\half z(\mu_j - i\nu_j)}$$ 
 $ C'_{k; ab c_1 c_2}$ is uniquely determined as the  coefficient of the monomial 
$$ \Pi_{j=1}^p (1 - e^{i z \alpha_j})^{-(a_j +1)} \cdot
\cdot \Pi_{j=1}^q\cdot
 (1 - e^{ z \lambda_j})^{-(b_j+1)} \cdot \Pi_{j=1}
(1 - e^{z (\mu_j + i\nu_j)})^{-(c_{1j} +1)}(1 - e^{z(\mu_j - i\nu_j)})^{-(c_{2j} +
1)}.$$ Since we are multiplying thru by a quantity independent of $a,b,c_1,c_2$,
it follows that $C'_{k; ab c_1 c_2}$ is uniquely determined as the coefficient the 
of the monomial
$$\Pi_{j=1}^p
 (1 - e^{i z \alpha_j})^{(k+1)-(a_j +1)}e^{\half i z \alpha_j} 
\cdot \Pi_{j=1}^q
 (1 - e^{ z \lambda_j})^{(k+1)-(b_j+1)}e^{\half z \lambda_j} \cdot$$
$$\cdot \Pi_{j=1}^{c}
(1 - e^{z (\mu_j + i\nu_j)})^{(k+1)-(c_{1j} +1)}(1 - e^{z(\mu_j -
i\nu_j)})^{(k+1)-(c_{2j} + 1)}e^{\half z (\mu_j +
i\nu_j)} e^{\half z(\mu_j - i\nu_j)}.$$
Expanding into an exponential polynomial, we find that $C'_{k; ab c_1 c_2}$ is
uniquely determined as the coefficient of the monomial
$$ \Pi_{j=1}^p  e^{i z ((k+1)-(a_j +1) + \half )\alpha_j}
\cdot \Pi_{j=1}^q e^{ z((k+1)-(b_j+1) + \half) \lambda_j} \cdot$$
$$\cdot \Pi_{j=1}^{c} e^{z ((k+1)-(c_{1j} +1)+ \half)(\mu_j + i\nu_j)} e^{z
((k+1)-(c_{2j}) +
\half)(\mu_j - i\nu_j)}).$$
Uniqueness follows from the fact that the vector $\beta + \half$ is a minimal element
of the set of exponent vectors occuring in  (6.5 k).
  Since $C'_{k; ab c_1 c_2} =C_{k; ab c_1
c_2}$ when $|a| + |b| + |c_1| + |c_2| = k+1$, we have determined 
$C_{k; ab c_1 c_2}.$

We then remove the $C_{k; a b c_1 c_2} 
(D_{\alpha}+ \half)^{a} (D_{\lambda}+ \half)^b
(D_{\mu + i
\nu}+ \half)^{c_1}  (  D_{\mu - i
\nu}+ \half)^{c_2}$  terms with $|a| + |b| + |c_1| + |c_2| = k+1$
 in (6.5 k). This leaves only terms with coefficients
$C_{k; a b c_1 c_2}$ with $|a| + |b| + |c_1| + |c_2| \leq k.$  Hence we can continue
the process of recovering coefficients until the end. \qed
\medskip

Let us now rewrite
  $$\sum_{\beta \in \Nb^{n}} c_{k;\beta} e^{ z \langle
\beta + (\half,\dots,\half), (\alpha,
\lambda, \mu, \nu) \rangle} \leqno(6.6a)$$
in the form
$$\sum_{j=1}^M a_{jk}  e^{z \omega_j}.\leqno(6.6b)$$
\medskip

\noindent{\bf (6.7) Lemma}~~~{\it The complex exponents $\omega_j$ in (6.6b),
together with $\pi$, are  independent over the rationals. Moreover, the
coefficients $c_{k:\beta}$ can be determined  from the
coefficients $a_{jk}$ }
\medskip

\noindent{\bf Proof}: The $\omega_j$'s are  rational linear combinations of
the exponents $\alpha_j, \lambda_j, \mu_j, \nu_j$, which  by assumption are
independent, with $\pi$, over the rationals. This independence also implies that the
exponents 
$\langle \beta + (\half,\dots,\half), (\alpha,\lambda, \mu, \nu) \rangle$
are all distinct.  Hence the coefficients in (6.6a)-(6.6b) are the same. 
\qed 
 \medskip

The proof of Theorem II is thus reduced to the following general statement about
exponential polynomials.
\medskip

\noindent{\bf (6.8) Lemma }{\it Suppose that the
 exponents $\omega_j$ of an exponential polynomial (6.6b) are independent (with
$\pi$) over the rationals.  Then the coefficients $a_{jk}$ of
  can be determined from the values of this
polynomial at $z = N \in \Nb.$}
\medskip

\noindent{\bf Proof}:  If not,  
there would exist a polynomial with
the given complex frequencies   which vanished at all integers $z=N$. 
But the  different terms $e^{z \omega_j}, e^{z \omega_k}$ have different
exponential growth rates along $z = N$ or $z = - N$ ($N \in \Nb$) unless $Re \omega_j
=  Re \omega_k.$  Let us write the large sum as a sum of smaller sums with a 
common $Re \omega.$  Each of the smaller sums
must separately vanish for $z \in \Nb.$  Multiply each one by the
relevant factor of $e^{- Re \omega}$. Each then turns  into an exponential
polynomial  with imaginary exponents,  which 
vanishes  for all $z = N.$  Since the frequencies are
independent(with
$\pi$) over ${\bf Q}$, each of these polynomials must vanish identically if it
vanishes at integral points.  Hence the coefficients $a_{jk}$ are uniquely 
determined by values of the large sum at integral points. \qed

\end{document}